\documentclass[11pt]{article}

\usepackage{amscd,amsmath}
\usepackage{mathrsfs}
\usepackage{amsfonts}
\usepackage{amssymb}
\usepackage{enumerate}
\usepackage{setspace}
\usepackage{color}

\usepackage[hyperindex]{hyperref}

\newtheorem{theorem}{Theorem}[section]
\newtheorem{lemma}{Lemma}[section]

\newtheorem{corollary}{Corollary}[section]

\newtheorem{remark}{Remark}[section]

\newcommand{\eqn}{\begin{eqnarray}}
\newcommand{\een}{\end{eqnarray}}

\numberwithin{equation}{section}

\begin{document}

\title{Analyticity of the Subcritical and Critical Quasi-Geostrophic equations in Besov Spaces}

\author{Hantaek Bae\footnote{Department of Mathematics, University of California, Davis. \ hantaek@math.ucdavis.edu}  \hspace{1cm} Animikh Biswas\footnote{Department of Mathematics and Statistics, University of Maryland, Baltimore County. \ abiswas@umbc.edu} \hspace{1cm} Eitan Tadmor\footnote{Department of Mathematics, Center for Scientific Computation and Mathematical Modeling, and Institute for Physical Sciences \& Technology, University of Maryland, College Park. \  tadmor@cscamm.umd.edu}}

\date{\empty}

\maketitle

\begin{abstract}
\noindent
We establish analyticity of the subcritical and critical quasi-geostrophic equations in critical Besov spaces. The main method is so-called Gevrey estimates, which is motivated by the work of Foias and Temam \cite{FoiasTemam}. We show that mild solutions are Gevrey regular, i.e. they satisfy the estimate 
\[
\sup_{t>0}\|e^{\alpha t^{1/\gamma}\Lambda_{1}}\theta(t)\|_{\cal L}<\infty
\]
for a suitably chosen $\alpha>0$ and a scaling invariant Besov space $\cal L$.
\end{abstract}


\section{Introduction and statement of main results} \label{sec:1}

In this paper, we study  analyticity properties of the dissipative quasi-geostrophic equations. These equations are derived from the more general quasi-geostrophic approximation for a rapidly rotating non-homogeneous fluid flow with small Rossby and Ekman numbers. For the case of special solutions with constant potential vorticity in the interior and constant buoyancy frequency (normalized to one), the general quasi-geostrophic equations reduce to the evolution equations for the temperature on the two-dimensional boundary. For  more details on physical properties  and derivation of these equations, see \cite{Constantin, ConstantinMajda, Pedlosky}.  Here we consider the initial value problem given by the following system of equations in two space dimensions: 
\begin{subequations} \label{eq:1.1}
\begin{align}
& \theta_{t} +v\cdot \nabla \theta +\kappa \Lambda^{\gamma}\theta=0,  \label{eq:1.1 a}\\
& v=\left(-\mathscr{R}_{2}\theta, \mathscr{R}_{1}\theta\right), \label{eq:1.1 b}\\
& \theta(t=0,x)=\theta_{0}(x), \label{eq:1.1 c}
\end{align}
\end{subequations}
where the scalar function $\theta$ is the potential temperature, $v$ is the fluid velocity,  and $\kappa$ is the viscosity coefficient. $\mathscr{R}=(\mathscr{R}_{1}, \mathscr{R}_{2})$ are the Riesz transforms whose symbols are given by $\displaystyle i\xi_{l}/|\xi|$ for $l=1,2$, and $\Lambda^{\gamma}$ is a Fourier multiplier whose symbol is given by $|\xi|^{\gamma}$. The cases $\gamma>1$, $\gamma=1$, and $\gamma <1$ are called respectively the subcritical, critical and supercritical quasi-geostrophic equations. In this paper, we study the subcritical and critical cases (i.e. $\gamma \ge 1$) and for simplicity we will set $\kappa=1$.  Formally, we can express a solution $\theta$ in the integral form:
\eqn \label{eq:1.2}
\theta(t)=e^{-t\Lambda^{\gamma}}\theta_{0} - \int^{t}_{0} \left[e^{-(t-s)\Lambda^{\gamma}} (v\cdot \nabla \theta)(s) \right]ds.
\een
Any solution satisfying this integral equation is called a \emph{mild solution}, which can be obtained by using a fixed point argument. An \emph{invariant space} for solving this integral equation corresponds to a scaling invariance property of the equation. This can be stated as follows. Assume that $\displaystyle (\theta, v) \in \big(\mathscr{L}, \left(-\mathscr{R}_{2}\mathscr{L}, \mathscr{R}_{1}\mathscr{L}\right)\big)$ solves (\ref{eq:1.1}) for a Banach space $\mathscr{L}$. Then, the same is true for rescaled functions:
\eqn \label{eq:1.3}
\theta_{\lambda}(t,x):=\lambda^{\gamma-1}\theta\left(\lambda^{\gamma} t,\lambda x \right), \quad v_{\lambda}(t,x):=\lambda^{\gamma-1}v\left(\lambda^{\gamma} t,\lambda x \right) \quad \forall \lambda>0.
\een
Under these scalings, $\mathscr{L}=L^{\frac{2}{\gamma-1}}$ ($\gamma \ge 1$), $\dot{H}^{2-\gamma}$, and $\dot{B}^{\frac{2}{p}+1-\gamma}_{p,q}$ are {\it  critical spaces} for initial data, i.e., the corresponding norms are invariants under these scalings. One can find various well-posedness results for small initial data in these critical spaces in \cite{Bae1, Bae2, Chae, Chen, Hmidi, Ju, Wu1, Wu2}. There are also several well-posedness results for large data;  in the Lebesgue space $L^{\frac{2}{\gamma-1}}$ for $\gamma>1$ (\cite{Carrillo}), in the energy space $H^{1}$ (\cite{Dong2}), and in Besov spaces (\cite{Abidi, Dong4, Wang}).

The goal of this paper is to show space analyticity of mild solutions and to provide explicit estimates for the analyticity radius as a function of time. In fluid-dynamics, the space analyticity radius has an important physical interpretation: at this length scale the viscous effects  and the (nonlinear) inertial effects  are roughly comparable.  Below this length scale, the viscous effects dominate the inertial effects and the Fourier spectrum decays exponentially  (\cite{Doelman, Foias, Henshaw1, Henshaw2}).  This fact can be used to show that the finite dimensional Galerkin approximations  converge exponentially fast in these cases (\cite{Doering}). Other applications of analyticity radius  occur in establishing sharp temporal decay rates of solutions in higher Sobolev norms  (\cite{Bae3, Oliver}), establishing geometric  regularity criteria for the Navier-Stokes equations, and in measuring the  spatial complexity of fluid flow (see \cite{Grujic, Kukavica1, Kukavica2}).

We follow here the Gevrey class approach pioneered by Foias and Temam (\cite{FoiasTemam}) for estimating space analyticity radius for the Navier-Stokes equations and which was subsequently used by many authors (\cite{Bae3, Biswas1, Biswas2, Cao, Ferrari}, and the references there in).  More precisely, we will show that mild solutions of (\ref{eq:1.1}) are { \it Gevrey regular}, i.e. they satisfy the estimate 
\[
\sup_{t>0}\left\|e^{\alpha t^{1/\gamma}\Lambda_{1}}\theta \right\|_{\cal L}<\infty,
\]
where $\cal L$ is a critical Besov space (which will be specified below) and $\alpha$ is a suitably chosen parameter (namely, $\alpha=1$ for $\gamma>1$ and $\alpha=\frac{1}{4}$ for $\gamma=1$). This approach enables  one to avoid cumbersome recursive estimation of higher order derivatives.   We emphasize that  the symbol of $\Lambda_{1}$  is given by the $l^{1}$ norm $|\xi|_{1}=|\xi_{1}|+|\xi_{2}|$ rather than the usual $\Lambda=\sqrt{-\Delta}$.

Before stating our main results, let us explain the main idea of the paper for the subcritical case ($\gamma>1$). The same approach will be used for the critical case. To this end, we define 
\eqn \label{eq:1.4}
\Theta=:e^{t^{1/\gamma}\Lambda_{1}}\theta, \quad V=(V_{1}, V_{2})=:e^{t^{1/\gamma}\Lambda_{1}} \left(v_{1}, v_{2}\right).
\een
Then, $V=(-\mathscr{R}_{2}\Theta, \mathscr{R}_{1}\Theta)$ and $(\Theta, V)$ satisfies the following equation:
\begin{equation} \label{eq:1.5}
 \begin{split}
 \Theta(t)&=e^{t^{1/\gamma}\Lambda_{1}-t\Lambda^{\gamma}}\theta_{0} - \int^{t}_{0}\left[ e^{t^{1/\gamma}\Lambda_{1}-(t-s)\Lambda^{\gamma} } \left(e^{-s^{1/\gamma}\Lambda_{1}}V \cdot \nabla e^{-s^{1/\gamma}\Lambda_{1}} \Theta \right)(s) \right]ds\\
 &=e^{t^{1/\gamma}\Lambda_{1}-t\Lambda^{\gamma}}\theta_{0} - \int^{t}_{0}\left[ e^{(t^{1/\gamma}-s^{1/\gamma})\Lambda_{1}}\nabla e^{-(t-s)\Lambda^{\gamma}} \cdot e^{s^{1/\gamma}\Lambda_{1}} \left(e^{-s^{1/\gamma}\Lambda_{1}}V  e^{-s^{1/\gamma}\Lambda_{1}} \Theta \right)(s) \right]ds,
 \end{split}
\end{equation}
where, we have used the fact that $V$ is divergence free and $\nabla$ commutes with any Fourier multiplier in order to obtain the last equality. Since $e^{t^{1/\gamma}|\xi|_{1}}$ is dominated by $e^{-t|\xi|^{\gamma}}$, the linear term, $e^{t^{1/\gamma}\Lambda_{1}-t\Lambda_{1}^{\gamma}}\theta_{0}$, closely resembles that of $\theta$,  namely, $e^{-t\Lambda_{1}^{\gamma}}\theta_{0}$. The estimates of the nonlinear term are similar to those of $\theta$ due to the boundedness property of the following bilinear operator $B_{s}$:
\eqn \label{eq:1.111}
B_{s}(f,g):=e^{s^{1/\gamma}\Lambda_{1}} \left(e^{-s^{1/\gamma}\Lambda_{1}} f  e^{-s^{1/\gamma}\Lambda_{1}} g\right).
\een
(For details, see Section \ref{sec:3.2}). 

As noticed from the above argument, the existence result of $\theta$ is crucial in establishing Gevrey regularity. Thus, in Section \ref{sec:3} (for $\gamma>1$) and Section \ref{sec:4} (for $\gamma=1$), we will first show the existence of a mild solution  and then proceed to explain how to modify the existence proof to obtain Gevrey regularity. We now present our existence/analyticity statements for $\gamma>1$ and $\gamma=1$ separately. For notational simplicity, we will suppress the dependence on $p$, $q$ and other relevant indices when we define norms below.

\subsection{Subcritical Case: $\gamma>1$}

\subsubsection{\bf Existence in Besov spaces}

We begin with the existence result in critical Besov spaces. Let us take initial data $\theta_{0}$ in $\dot{B}^{\frac{2}{p}+1-\gamma}_{p,q} \cap L^{p_{1}}$, with 
\eqn  \label{paramres}
1 \leq q \leq \infty,  \quad p_{1} = \frac{2}{\gamma-1} \quad \mbox{and}\  p < \frac{2}{\gamma-1}.
\een 
Since the dissipation rate $\gamma>1$ dominates the derivative in the advection term, we can follow Weissler's idea (\cite{Weissler}) to define a function space $E_{\gamma,T}=K_{\gamma,T}\cap G_{\gamma,T}$, with 
\begin{equation*}
 \begin{split}
 & \|\theta\|_{K_{\gamma, T}} :=\sup_{0<t\leq T}\left[ \|\theta(t)\|_{\dot{B}^{\frac{2}{p}+1-\gamma}_{p,q}}  +t^{\alpha/ \gamma} \|\theta(t)\|_{\dot{B}^{\frac{2}{p}+1-\gamma+\alpha}_{p,q}}\right],\\
 &\|\theta\|_{G_{\gamma, T}} :=\sup_{0<t<T}\left[\|\theta(t)\|_{L^{p_{1}}} +t^{\beta}\|\theta(t)\|_{L^{r}}\right].
 \end{split}
\end{equation*}
We note that the time weights $t^{\alpha/\gamma}$ and $t^{\beta}$ indicate the gain  in regularity and integrability of the solution through the dissipation operator $\Lambda^{\gamma}$:
\[
|\xi|^{a}e^{-t|\xi|^{\gamma}} \lesssim t^{-\frac{a}{\gamma}}.
\]
The parameters $r, \alpha $ and $\beta$ are chosen to satisfy
\eqn  \label{alphachoice}
1<r<\infty, \quad \frac{2}{r} < \min\left\{ \gamma -1, \frac{2}{p} - (\gamma-1)\right\}, \quad \alpha < 1 + \frac{2}{r}
\quad\mbox{and}\ \beta=1-1/\gamma-\frac{2}{r\gamma}.
\een 
It is possible to choose $r$ and $\alpha$ as in (\ref{alphachoice}) due to the fact that $\gamma >1$ and the restriction on the parameters given in
(\ref{paramres}). With the above mentioned choices, we must have the relations
\eqn  \label{eq:1.6}
(i) \hspace{0.2cm} 0<\frac{2}{r}< \gamma-1, \quad (ii) \hspace{0.2cm} \frac{2}{p}+1-\gamma-\frac{2}{r}>0, \quad 
(iii) \hspace{0.2cm} 0<\beta <1-\frac{\alpha}{\gamma}.
\een
(Since the condition (\ref{eq:1.6}) (or equivalently, (\ref{paramres})) will be used to prove the existence and analyticity results stated below (Theorems \ref{thm:1.1}-\ref{thm:1.3}), we will assume it for these results.) Condition (i) implies $\beta>0$, and conditions (ii) and (iii) in (\ref{eq:1.6}), which lead to the restriction $p < \frac{2}{\gamma-1}$ in (\ref{paramres}), are  necessary to estimate the quadratic term $v\theta$ in the high frequency part, while, as in \cite{Carrillo}, the condition $p_{1}= \frac{2}{\gamma-1}$ is required to  control the low frequency part.   To prove the existence of a solution in $E_{\gamma, T}$, we first need the existence of a solution in $G_{\gamma,T}$.

\begin{theorem} \cite{Carrillo} \label{thm:carrillo}
For any initial data $\theta_{0} \in L^{p_{1}}$ with $p_{1}= \frac{2}{\gamma-1}$, there exists a global-in-time solution $\theta$ in $G_{\gamma,T}$. Moreover, for any time interval $[t, t+\tau]$,
\[
\lim_{\tau\rightarrow 0} \left[\tau^{\beta}\|\theta(t+\tau)\|_{L^{r}}\right] =0 \quad \text{uniformly in} \ t.
\]
\end{theorem}

We now state our first existence result of this paper.  

\begin{theorem} [Existence] \label{thm:1.1}
For any initial data $\theta_{0} \in \dot{B}^{\frac{2}{p}+1-\gamma}_{p,q}\cap L^{p_{1}}$, there exists a global-in-time solution $\theta \in K_{\gamma,T}$ such that
\eqn\label{eq:1.7}
\left\|\theta\right\|_{K_{\gamma, T}} \lesssim 2^{T} \left\|\theta_{0}\right\|_{\dot{B}^{\frac{2}{p}+1-\gamma}_{p,q}}. 
\een
\end{theorem}

As one can see in Section \ref{sec:3}, the main step is to show 
\[
\left\|\theta \right\|_{K_{\gamma, T}} \lesssim \left\|\theta_{0} \right\|_{\dot{B}^{\frac{2}{p}+1-\gamma}_{p,q}} + \sup_{0< \tau \leq T}\left[\tau^{\beta}\|\theta(\tau)\|_{L^{r}}\right] \left\|\theta \right\|_{K_{\gamma, T}}.
\]
Therefore, Theorem \ref{thm:carrillo} is crucial to prove the following result. We note that the size of initial data is arbitrary.

\subsubsection{\bf Analyticity in Besov spaces}
 As we use the existence of $\theta$ in $G_{\gamma,T}$ to show the existence proof of $\theta$ in $K_{\gamma,T}$, the proof of analyticity of the solution of Theorem \ref{thm:1.1} exactly consists of two parts. We first show the existence of $\Theta$ in $G_{\gamma,T}$ (Theorem \ref{thm:1.2}) which provides the smallness of $\tau^{\beta}\|\Theta(t+\tau)\|_{L^{q}}$ as $\tau \rightarrow 0$. Then, we show the existence of $\Theta$ in $K_{\gamma, T}$ (Theorem \ref{thm:1.3}) by obtaining 
\[
\left\|\Theta \right\|_{K_{\gamma, T}} \lesssim \left\|\theta_{0} \right\|_{\dot{B}^{\frac{2}{p}+1-\gamma}_{p,q}} + \sup_{0< \tau \leq T}\left[\tau^{\beta}\|\Theta(\tau)\|_{L^{r}}\right] \|\Theta\|_{K_{\gamma, T}}.
\]
We note that the proof of these results is based on the boundedness property of the bilinear operator $B_{s}$ defined in (\ref{eq:1.111}). 

\begin{theorem} [Existence]\label{thm:1.2}
For any initial data $\theta_{0} \in L^{p_{1}}$ with $p_{1}= \frac{2}{\gamma-1}$, there exists a global-in-time solution $\Theta(t)\in G_{\gamma,T}$. Moreover, for any time interval $[t, t+\tau]$,
\[
\lim_{\tau\rightarrow 0} \left[\tau^{\beta}\|\Theta(t+\tau)\|_{L^{r}}\right] =0 \quad \text{uniformly in} \ t.
\]
\end{theorem}

\begin{theorem} [Analyticity] \label{thm:1.3}
For any initial data $\theta_{0} \in \dot{B}^{\frac{2}{p}+1-\gamma}_{p,q}\cap L^{p_{1}}$, there exists a global-in-time solution $\Theta \in K_{\gamma, T}$ such that
\eqn \label{eq:1.9}
\left\|\Theta \right\|_{K_{\gamma, T}} \lesssim 2^{T} \left\|\theta_{0} \right\|_{\dot{B}^{\frac{2}{p}+1-\gamma}_{p,q}}. 
\een
\end{theorem}

\begin{remark}\label{rem:1.1}
Theorem \ref{thm:1.2} provides a new proof of the results in \cite{Dong3}, where they prove analyticity of solutions of the subcritical quasi-geostrophic equations by recursive estimation of higher order derivatives which involves a sophisticated combinatorial argument. By contrast, we provide analyticity directly from the mild solution setting and the proof is significantly reduced. Moreover, our method is more flexible and can be applied to other function spaces as Theorem \ref{thm:1.3}.
\end{remark}

\subsection{Critical case $\gamma=1$}

\subsubsection{\bf Existence in critical spaces}
We now consider the critical case. Unlike to  the subcritical case where the dissipation is strong enough to overcome the derivative in the advection term,  main ingredients for proving the global well-posedness for the critical case are the commutator estimates and the modulus of continuity argument (\cite{Abidi, Dong2, Dong4, Wang}.) However, we cannot apply Gevrey regularity approach to these work because  the exponential operator $\displaystyle e^{\alpha t \Lambda_{1}}$ does not commute with the advection term $v\cdot \nabla$. Thus, we need a result which does not rely on commutator estimates, which is the case of \cite{Bae2}. We will use this result to establish Gevrey regularity.

\begin{theorem} [Existence \cite{Bae2}] \label{thm:1.4}
There exists a constant $\epsilon_{0}>0$ such that for all initial data $\theta_{0} \in L^{\infty} \cap \dot{H}^{1} \cap \dot{B}^{\frac{2}{p}}_{p,q}$ and $v_{0} \in L^{\infty}$, with
\[
\left\|\theta_{0}\right\|_{L^{\infty}}+\left\|v_{0}\right\|_{L^{\infty}} +\|\theta_{0}\|_{\dot{H}^{1}}<\epsilon_{0},
\]
there exists a global-in-time solution $\theta \in L^{\infty}_{t}(L^{\infty} \cap \dot{H}^{1}) \cap L^{2}_{t}\dot{H}^{\frac{3}{2}} \cap E_{1}$, where
\eqn \label{eq:1.10}
\left\|\theta\right\|_{E_{1}}:=\left\|\theta\right\|_{\tilde{L}^{\infty}_{t}\dot{B}^{\frac{2}{p}}_{p,q}}+ \left\|\theta\right\|_{\tilde{L}^{1}_{t} \dot{B}^{\frac{2}{p}+1}_{p,q}}.
\een
\end{theorem}

The proof of Theorem \ref{thm:1.4} consists of two parts. We first show the existence of small solutions in $L^{\infty}_{t}(L^{\infty} \cap \dot{H}^{1}) \cap L^{2}_{t}\dot{H}^{\frac{3}{2}}$ with small initial data in $L^{\infty}\cap \dot{H}^{1}$. Since the Riesz transforms $\mathscr{R}$ defining $v$ do not map $L^{\infty}$ to $L^{\infty}$, we thus need to assume $v_{0}\in L^{\infty}$ as well. To estimate $v$ in $L^{\infty}$, we will use the Oseen kernel representation. (The definition of the Oseen kernel and its property will be presented in Section \ref{sec:4}). Then, using this smallness condition in $L^{\infty} \cap \dot{H}^{1}$, we prove the existence of  a solution in $E_{1}$. We note that depending on the ranges of $p$ and $q$, we have different restrictions on the size of $\theta_{0}$ in $\dot{B}^{\frac{2}{p}}_{p,q}$. If $\dot{B}^{\frac{2}{p}}_{p,q} \subset \dot{H}^{1}$, the size of $\|\theta_{0}\|_{\dot{B}^{\frac{2}{p}}_{p,q}}$ can be arbitrary, while if $\dot{H}^{1} \subset \dot{B}^{\frac{2}{p}}_{p,q}$, the size of $\|\theta_{0}\|_{\dot{B}^{\frac{2}{p}}_{p,q}}$ need to be small.

\subsubsection{\bf Analyticity in critical spaces} 
To show that solutions in Theorem \ref{thm:1.4} are analytic, we need to define the Gevrey operator more carefully because $\Lambda$ and $\Lambda_{1}$ are equivalent as a Fourier multiplier. Let
\[
\Theta=e^{\frac{1}{4}t\Lambda_{1}}\theta(t), \quad V=e^{\frac{1}{4}t\Lambda_{1}}v(t).
\]
Then, $\Theta$ and $V$ satisfy the following equation;
\eqn \label{eq:1.11}
\Theta(t)=e^{\frac{1}{4}t\Lambda_{1}-t\Lambda}\theta_{0}-\int^{t}_{0}\left[ e^{\frac{1}{4}t\Lambda_{1}-(t-s)\Lambda} \left(e^{-\frac{1}{4}s\Lambda_{1}}V \cdot \nabla e^{-\frac{1}{4}s\Lambda_{1}} \Theta\right)(s) \right]ds.
\een
As we will see in Section \ref{sec:4}, the fact $\frac{1}{4}|\xi|_{1} <\frac{1}{2}|\xi|$ allows that (\ref{eq:1.11}) is equivalent to the integral equation of $\theta$. Therefore, we can prove the following result along the lines of the proof of Theorem \ref{thm:1.3} and Theorem \ref{thm:1.4}.

\begin{theorem} [Analyticity] \label{thm:1.5}
There exists a constant $\epsilon_{0}>0$ such that for all initial data $\theta_{0} \in L^{\infty} \cap \dot{H}^{1} \cap \dot{B}^{\frac{2}{p}}_{p,q}$ and $v_{0} \in L^{\infty}$, with
\[
\left\|\theta_{0}\right\|_{L^{\infty}}+\left\|v_{0}\right\|_{L^{\infty}} +\left\|\theta_{0}\right\|_{\dot{H}^{1}}<\epsilon_{0},
\]
there exists a global-in-time solution $\Theta \in L^{\infty}_{t}(L^{\infty} \cap \dot{H}^{1}) \cap L^{2}_{t}\dot{H}^{\frac{3}{2}} \cap E_{1}$.
\end{theorem}

Theorem \ref{thm:1.5} in fact provides the decay estimates of the solution in Theorem \ref{thm:1.4}. 

\begin{corollary}
The solution in Theorem \ref{thm:1.4} satisfies the following decay rate of the solution;
\[
\left\|\nabla^{k} \theta(t)\right\|_{L^{\infty}\cap \dot{H}^{1}\cap \dot{B}^{\frac{2}{p}}_{p,q}} \leq  C_{k} t^{-k}, \quad \forall k>0.
\]
\end{corollary}

The proof  is rather simple so we provide it here. Since $L^{\infty}\cap \dot{H}^{1}\cap \dot{B}^{\frac{2}{p}}_{p,q}$ is translation-invariant and the operator $\nabla^{k} e^{-t \Lambda_{1}} $ is $L^{1}$-bounded,
\begin{equation}
 \begin{split}
\left\|\nabla^{k} \theta(t) \right\|_{L^{\infty}\cap \dot{H}^{1}\cap \dot{B}^{\frac{2}{p}}_{p,q}}  &=\left\|\nabla^{k} e^{-t^{1/\gamma}\Lambda_{1}} e^{t^{1/\gamma}\Lambda_{1}} \theta(t) \right\|_{L^{\infty}\cap \dot{H}^{1}\cap \dot{B}^{\frac{2}{p}}_{p,q}}\\
 & \leq \left\|\nabla^{k} e^{-t^{1/\gamma}\Lambda_{1}} \right\|_{L^{1}} \left\|\Theta(t)\right\|_{L^{\infty}\cap \dot{H}^{1}\cap \dot{B}^{\frac{2}{p}}_{p,q}}  \leq C \left\|\nabla^{k} e^{-t^{1/\gamma}\Lambda_{1}}\right\|_{L^{1}}
 \end{split}
\end{equation}
where $C=\left\|\Theta(t)\right\|_{L^{\infty}\cap \dot{H}^{1}\cap \dot{B}^{\frac{2}{p}}_{p,q}} \lesssim \left\|\theta_{0}\right\|_{L^{\infty}\cap \dot{H}^{1}\cap \dot{B}^{\frac{2}{p}}_{p,q}}$. Then, by rescaling $\xi$ as $\xi\mapsto t\xi$,
\[
\left\|\nabla^{k} e^{-t \Lambda_{1}} \right\|_{L^{1}}  \lesssim C_{k} t^{-k}, \quad C_{k}=\left\|\nabla^{k} e^{-\Lambda_{1}}\right\|_{L^{1}}.
\]

\begin{remark} 
(1) In \cite{Dong4}, the global well-posedness of the critical quasi-geostrophic equations is proved for large data in $\dot{B}^{\frac{2}{p}}_{p,q}$. However, analyticity of the solution was not proven; also the decay estimates provided there were for short times only. \\
(2) In \cite{Constantin}, they proved that for initial data $\theta_{0}\in H^{2}(\mathbb{T}^{2})$ with the smallness condition to $\|\theta_{0}\|_{L^{\infty}}$ there exists $t_{0}>0$ such that the solutions of the critical quasi-geostrophic equations are analytic for $t>t_{0}$. Along the proof in \cite{Constantin}, Dong proved in \cite{Dong1} that there exists a time $t_{0}>0$ such that the solution of the critical quasi-geostrophic equations is analytic for all $t\ge t_{0}$ for large data in $\dot{H}^{1}(\mathbb{T}^{2})$. Although we need smallness condition in $L^{\infty}\cap \dot{H}^{1}$, the solution of Theorem \ref{thm:1.4} becomes analytic instantaneously and it covers more spaces for initial data due to the range of $p$ and $q$.
\end{remark}

\section{Notations: Littlewood-Paley decomposition and paraproduct}\label{sec:2}
We begin with some notation. 
\begin{enumerate}[]
\item (1) Let $X$ be a  Banach space. $L^{p}(0,T;X)$ denotes the Banach space of Bochner measurable functions $f$ from $(0,T)$ to $X$ endowed with either the norm $\left(\int ^{T}_{0} \|f(\cdot, t)\|^{p}_{X}dt\right)^{\frac{1}{p}}$ for $1\leq p<\infty$ or $\displaystyle\sup_{0\leq t\leq T} \|f(\cdot,t)\|_{X}$ for $p=\infty$. For $T=\infty$, we use the notation $L^{p}_{t}X$ instead of $L^{p}(0,\infty;X)$. 
\item (2) For a sequence $\{a_{j}\}_{j\in\mathbb{Z}}$, $\{a_{j}\}_{l^{q}}:=\big(\displaystyle\sum_{j\in\mathbb{Z}} |a_{j}|^{q}\big)^{\frac{1}{q}}$, with the usual modification for $q=\infty$.
\item (3) $A \lesssim B $ means there exists a constant $C>0$ such that $A \leq C B $. Similarly, $A \sim B$ means that there exist two constants $C_1,C_2 >0$ such that $A \le C_1B$ and $B\le C_2A$.
\end{enumerate}

We next provide notation and definitions in the Littlewood-Paley theory. We take a couple of smooth functions $(\chi, \phi)$ supported on $\{\xi ; |\xi|\leq 1\}$ with values in $[0,1]$ such that for all $\xi \in \mathbb{R}^{d}$,
\[
\chi(\xi)+ \sum^{\infty}_{j=0}\psi(2^{-j}\xi)=1, \quad \psi(\xi)=\phi \left(\xi/2\right)-\phi(\xi).
\]
We set $\psi(2^{-j}\xi)=\psi_{j}(\xi)$ and define the dyadic blocks and lower frequency cut-off functions by 
\eqn\label{eq:2.1}
\Delta_{j}u=2^{jd} \int_{R^{d}} h(2^{j}y)u(x-y)dy, \quad S_{j}u=2^{jd} \int_{R^{d}} \tilde{h}(2^{j}y)u(x-y)dy,
\een
with $h=\mathscr{F}^{-1}\psi$ and $\tilde{h}=\mathscr{F}^{-1}\chi$. We define the {\it homogeneous Littlewood-Paley decomposition} by 
\eqn \label{eq:2.2}
u=\sum_{j\in \mathbb{Z}} \Delta_{j} u \quad \text{in} \quad \mathscr{S}^{'}_{h},
\een
where $\mathscr{S}^{'}_{h}$ is the space of tempered distributions $u$ such that $ \displaystyle \lim_{j\rightarrow -\infty}S_{j}u=0$ in $\mathscr{S}^{'}$. Using this decomposition, we define {\it homogeneous Besov spaces} as follows:
\[
\dot{B}^{s}_{p,q}=\Big\{f\in \mathscr{S}^{'}_{h} \ ; \ \|f\|^{q}_{\dot{B}^{s}_{p,q}}:=\sum_{j\in \mathbb{Z}}2^{jsq}\left\|\Delta_{j}f\right\|^{q}_{L^{p}}<\infty \Big\}.
\]
We also define time dependent {\it homogeneous Besov spaces} ;
\begin{equation} \label{eq:2.3}
 \begin{split}
 & L^{r}(0, T;\dot{B}^{s}_{p,q}) =\left\{f\in \mathscr{S}^{'}_{h} \ ; \ \|f\|_{L^{r}(0, T;\dot{B}^{s}_{p,q})}:=\left\|\left(\sum_{j\in \mathbb{Z}} 2^{jsq}\|\Delta_{j}f\|^{q}_{L^p} \right)^{\frac{1}{q}} \right\|_{L^{r}(0,T)}<\infty \right\}, \\
 & \tilde{L}^{r}(0,T; \dot{B}^{s}_{p,q}) =\left\{f\in \mathscr{S}^{'}_{h} \ ; \ \|f\|_{\tilde{L}^{r}(0,T; \dot{B}^{s}_{p,q})}: = \left(\sum_{j\in \mathbb{Z}} 2^{jsq}\left\|\Delta_{j}f\right\|^{q}_{L^{r}(0,T; L^p)}\right)^{\frac{1}{q}}<\infty \right\}
 \end{split}
\end{equation}
with the usual modification for $q=\infty$.

The concept of {\it  paraproduct} enables to deal with the interaction of two functions in terms of low or high frequency parts, \cite{Chemin}. For two tempered distributions $f$ and $g$,
\begin{equation} \label{eq:2.5}
 \begin{split}
 & fg=T_{f}g + T_{g}f + R(f,g),\\
 & T_{f}g= \sum_{i\leq j-2}\Delta_{i}f\Delta_{j}g= \sum_{j\in\mathbb{Z}}S_{j-1}f\Delta_{j}g, \quad R(f,g)= \sum_{|j-j^{'}|\leq 1}\Delta_{j}f \Delta_{j^{'}}g.
 \end{split}
\end{equation}
Then, up to finitely many terms,
\eqn \label{eq:2.6}
\Delta_{j}(T_{f}g)=S_{j-1}f \Delta_{j}g, \quad \Delta_{j}R(f,g)=\sum_{k\ge j-2}\Delta_{k}f\Delta_{k}g.
\een

We finally recall a few inequalities which will be used in the sequel.
\begin{lemma}
(1) Bernstein's inequality \cite{Chemin}: For $1\leq p\leq q \leq \infty$ and $k\in \mathbb{N}$,
\eqn \label{eq:2.7}
\sup_{|\alpha|=k} \left\|\partial^{\alpha}\Delta_{j}f \right\|_{L^{p}} \simeq 2^{jk} \left\|\Delta_{j}f \right\|_{L^{p}}, \quad \left\|\Delta_{j}f \right\|_{L^{q}} \lesssim 2^{jd(\frac{1}{p}-\frac{1}{q})} \left\|\Delta_{j}f \right\|_{L^{p}}.
\een
(2) Localization of the fractional heat kernel \cite{Hmidi}: 
\eqn \label{eq:2.8}
\left\|e^{-t\Lambda_{1}^{\gamma}} \Delta_{j}f \right\|_{L^{p}} \lesssim e^{-t2^{\gamma j}} \left\|\Delta_{j}f\right\|_{L^{p}}.
\een
\end{lemma}

\section{Proof of Theorem \ref{thm:1.1}, Theorem \ref{thm:1.2} and Theorem \ref{thm:1.3}}\label{sec:3}
In this section, we prove the existence and analyticity of the subcritical quasi-geostrophic equations with large initial data in $\dot{B}^{\frac{2}{p}+1-\gamma}_{p,q}\cap L^{p_{1}}$. As we already mentioned in the introduction, we need Theorem \ref{thm:carrillo} to show the existence of a solution in $K_{\gamma,T}$. We recall the definition of the norm of $G_{\gamma,T}$ and $K_{\gamma,T}$;
\begin{equation} \label{eq:3.1}
 \begin{split}
 & \|\theta\|_{G_{\gamma,T}} =\sup_{0<t<T}\left[\|\theta(t)\|_{L^{p_{1}}} +t^{\beta}\|\theta(t)\|_{L^{r}}\right], \quad \beta=1-1/\gamma-\frac{2}{r\gamma}>0,\\
 & \|\theta\|_{K_{\gamma,T}} =\sup_{0<t\leq T}\left[ \|\theta(t)\|_{\dot{B}^{\frac{2}{p}+1-\gamma}_{p,q}}  +t^{\alpha/\gamma} \|\theta(t)\|_{\dot{B}^{\frac{2}{p}+1-\gamma+\alpha}_{p,q}}\right].
 \end{split}
\end{equation}
In addition, we will use the following lemma repeatedly in the proof of Theorem \ref{thm:1.1}.
\begin{lemma} \label{lem:3.1}
For any $0<a<1$ and $0<b<1$,
\[
\int^{t}_{0} \Big[ (t-s)^{-a} s^{-b} \Big] ds \lesssim t^{1-a-b}.
\]
\end{lemma}

\subsection{Proof of Theorem \ref{thm:1.1}} \label{sec:3.1}
We first show that a priori estimate
\eqn \label{eq:3.4}
\|\theta\|_{K_{\gamma, T}} \lesssim \|\theta_{0}\|_{\dot{B}^{\frac{2}{p}+1-\gamma}_{p,q}} + \sup_{0< \tau \leq T}\left[\tau^{\beta}\|\theta(\tau)\|_{L^{r}}\right] \|\theta\|_{K_{\gamma, T}}
\een
is enough to prove Theorem \ref{thm:1.1}. The proof comes as follows. We take $T_{1}$ such that on $[0,T_{1}]$
\[
\tau^{\beta}\|\theta(\tau)\|_{L^{r}} \leq \frac{1}{2} \ \Longrightarrow \|\theta\|_{K_{\gamma, t}}\lesssim  2 \|\theta_{0}\|_{\dot{B}^{\frac{2}{p}+1-\gamma}_{p,q}}.
\]
By Theorem \ref{thm:carrillo}, we can take the next step by taking $T_{2}=2T_{1}$. Then, (\ref{eq:3.4}) implies that 
\[
\|\theta(t)\|_{\dot{B}^{\frac{2}{p}+1-\gamma}_{p,q}} \lesssim  4 \|\theta_{0}\|_{\dot{B}^{\frac{2}{p}+1-\gamma}_{p,q}}
\]
on $[T_{1}, 2T_{1}]$. Inductively, we can obtain that for $t\in [nT_{1}, (n+1)T_{1}]$
\[
\|\theta\|_{K_{\gamma, t}} \lesssim 2^{(n+1)} \|\theta_{0}\|_{\dot{B}^{\frac{2}{p}+1-\gamma}_{p,q}}
\]
which leads to (\ref{eq:1.7}). To prove (\ref{eq:3.4}), we express a solution $\theta$ in the integral form:
\eqn \label{eq:3.5}
\theta(t)=e^{-t\Lambda^{\gamma}}\theta_{0} - \int^{t}_{0} \left[\nabla e^{-(t-s)\Lambda^{\gamma}} \cdot (v \theta)(s) \right]ds := e^{-t\Lambda^{\gamma}}\theta_{0}- \mathscr{B}(v,\theta).
\een
We note that the definition of the space $K_{\gamma,T}$ is based on the linear behavior of the solution. So, it  is easy to show that 
\eqn \label{linear}
\left\|e^{-t\Lambda^{\gamma}}\theta_{0}\right\|_{K_{\gamma,T}} \lesssim \|\theta_{0}\|_{\dot{B}^{\frac{2}{p}+1-\gamma}_{p,q}}.
\een
Therefore, it is enough to estimate $\mathscr{B}(v,\theta)$ in $K_{\gamma,T}$. Let us decompose $v \theta$ as a paraproduct: $ v \theta =T_{v}\theta +T_{\theta}v +R(v,\theta)$. Then,
\begin{equation*}
 \begin{split}
 \mathscr{B}(v,\theta)&=\int^{t}_{0} \left[\nabla e^{-(t-s)\Lambda^{\gamma}} \cdot (v \theta)(s)\right] ds \\
 &= \int^{t}_{0} \left[\nabla e^{-(t-s)\Lambda^{\gamma}} \cdot \left(T_{v}\theta +T_{\theta} v +R(v,\theta)\right) \right]ds\\
 &:=\mathscr{B}_{1}(v,\theta)+\mathscr{B}_{2}(v,\theta)+\mathscr{B}_{3}(v,\theta).
 \end{split}
\end{equation*}
In the sequel, we will treat $v$ as $\theta$ in the estimations of $\mathscr{B}_{i}(v,\theta)$'s in $K_{\gamma, T}$ because $v=(-\mathscr{R}_{2}\theta, \mathscr{R}_{1}\theta)$ and the Riesz transforms are bounded in $L^{p}$ for all $p\in(1, \infty)$.

\subsubsection{\bf Estimation of $\mathscr{B}_{1}(v,\theta)$ and $\mathscr{B}_{2}(v,\theta)$} 
We take the dyadic operator $\Delta_{j}$ to $\mathscr{B}_{1}(v,\theta)$ and take the $L^{p}$ norm. By (\ref{eq:2.6}) and (\ref{eq:2.8}), we have
\eqn \label{eq:3.6}
\left\|\Delta_{j}\mathscr{B}_{1}(v,\theta)(t) \right\|_{L^{p}} \lesssim \int^{t}_{0} \left[2^{j} e^{-(t-s)2^{\gamma j}} \left\|S_{j-1}\theta(s) \right\|_{L^{\infty}} \left\|\Delta_{j} \theta(s) \right\|_{L^{p}} \right]ds.
\een
We first estimate $\displaystyle \left\|S_{j-1}\theta\right\|_{L^{\infty}}$: 
\eqn \label{eq:3.7}
\left\|S_{j-1}\theta(s) \right\|_{L^{\infty}} \lesssim \sum^{j-1}_{k=-\infty} 2^{\frac{2k}{r}} \left\|\Delta_{k} \theta(s) \right\|_{L^{r}} \lesssim 2^{\frac{2j}{r}}s^{-\beta}s^{\beta} \left\|\theta\right\|_{L^{r}},
\een
where we use (\ref{eq:2.7}) for the first inequality. Thus, the right-hand side of (\ref{eq:3.6}) can be bounded by
\begin{equation} \label{eq:3.8}
 \begin{split}
 \left\|\Delta_{j}\mathscr{B}_{1}(v,\theta)(t) \right\|_{L^{p}} &\lesssim \int^{t}_{0} \left[2^{j} e^{-(t-s)2^{\gamma j}} \left\|S_{j-1}\theta(s)\right\|_{L^{\infty}} \left\|\Delta_{j} \theta(s)\right\|_{L^{p}} \right]ds\\
 & \lesssim \sup_{0< \tau \leq t}\left[\tau^{\beta}\|\theta(\tau)\|_{L^{r}}\right] \int^{t}_{0} \left[2^{j(1+\frac{2}{r})} e^{-(t-s)2^{\gamma j}} s^{-\beta} \|\Delta_{j}\theta(s)\|_{L^{p}}\right]ds.
 \end{split}
\end{equation}
Since 
\[
2^{j(1+\frac{2}{r})} e^{-(t-s)2^{\gamma j}} \lesssim (t-s)^{-1/\gamma(1+\frac{2}{r})},
\]
we have
\eqn \label{eq:3.9}
\left\|\Delta_{j}\mathscr{B}_{1}(v,\theta)(t) \right\|_{L^{p}} \lesssim \sup_{0< \tau \leq t}\left[\tau^{\beta}\|\theta(\tau)\|_{L^{r}}\right] \int^{t}_{0} \left[(t-s)^{-1/\gamma(1+\frac{2}{r})}s^{-\beta} \left\|\Delta_{j}\theta(s)\right\|_{L^{p}}\right]ds.
\een
We multiply (\ref{eq:3.9}) by $2^{j \left(\frac{2}{p}+1-\gamma\right)}$ and take the $l^{q}$ norm. Then,
\begin{equation}\label{eq:3.10}
 \begin{split}
& \left\|\mathscr{B}_{1}(v,\theta)(t) \right\|_{\dot{B}^{\frac{2}{p}+1-\gamma}_{p,q}}  \\
&\lesssim \sup_{0< \tau \leq t}\left[\tau^{\beta}\|\theta(\tau)\|_{L^{r}}\right] \sup_{0<\tau<t}\left[\|\theta(\tau)\|_{\dot{B}^{\frac{2}{p}+1-\gamma}_{p,q}}\right] \int^{t}_{0} \left[(t-s)^{-1/\gamma(1+\frac{2}{r})}s^{-\beta} \right]ds  \\
 & \lesssim \sup_{0< \tau \leq t}\left[\tau^{\beta}\|\theta(\tau)\|_{L^{r}}\right] \left\|\theta\right\|_{K_{\gamma,t}} \int^{t}_{0} \left[(t-s)^{-1/\gamma(1+\frac{2}{r})}s^{-\beta} \right]ds \\
 &\lesssim \sup_{0< \tau \leq t}\left[\tau^{\beta}\|\theta(\tau)\|_{L^{r}}\right]  \left\|\theta\right\|_{K_{\gamma,t}},
 \end{split}
\end{equation}
where we use $\beta=1-1/\gamma-\frac{2}{\gamma r}$ and Lemma \ref{lem:3.1}  to bound the time integration. 

We next estimate $\|\mathscr{B}_{1}(v,\theta)\|_{\dot{B}^{\frac{2}{p}+1-\gamma+\alpha}_{p,q}}$. From (\ref{eq:3.9}),
\begin{equation}\label{eq:3.11}
 \begin{split}
 &2^{j \left(\frac{2}{p}+1-\gamma+\alpha\right)} \left\|\Delta_{j}\mathscr{B}_{1}(v,\theta)(t)\right\|_{L^{p}} \\
 & \lesssim \sup_{0< \tau \leq t}\left[\tau^{\beta}\|\theta(\tau)\|_{L^{r}}\right] \int^{t}_{0} \left[(t-s)^{-1/\gamma(1+\frac{2}{r})}s^{-\beta} 2^{j (\frac{2}{p}+1-\gamma+\alpha)} \left\|\Delta_{j}\theta(s) \right\|_{L^{p}}\right]ds\\
 &=\sup_{0< \tau \leq t}\left[\tau^{\beta}\|\theta(\tau)\|_{L^{r}}\right] \int^{t}_{0} \left[(t-s)^{-1/\gamma(1+\frac{2}{r})}s^{-\beta} s^{-\frac{\alpha}{\gamma}} s^{\frac{\alpha}{\gamma}}2^{j (\frac{2}{p}+1-\gamma+\alpha)} \left\|\Delta_{j}\theta(s) \right\|_{L^{p}}\right]ds.
 \end{split}
\end{equation}
By taking the $l^{q}$ norm, we have
\begin{equation}\label{eq:3.12}
 \begin{split}
 & \left\|\mathscr{B}_{1}(v,\theta)(t) \right\|_{\dot{B}^{\frac{2}{p}+1-\gamma+\alpha}_{p,q}} \\
 & \lesssim \sup_{0< \tau \leq t}\left[\tau^{\beta}\|\theta(\tau)\|_{L^{r}}\right] \sup_{0<\tau<t}\left[\tau^{\frac{\alpha}{\gamma}} \left\|\theta(\tau)\right\|_{\dot{B}^{\frac{2}{p}+1-\gamma+\alpha}_{p,q}}\right] \int^{t}_{0} \left[(t-s)^{-1/\gamma(1+\frac{2}{q})} s^{-\beta} s^{-\frac{\alpha}{\gamma}} \right]ds \\
 & \lesssim \sup_{0< \tau \leq t}\left[\tau^{\beta}\|\theta(\tau)\|_{L^{r}}\right] \left\|\theta\right\|_{K_{\gamma,t}}\int^{t}_{0} \left[(t-s)^{-1/\gamma(1+\frac{2}{q})} s^{-\beta} s^{-\frac{\alpha}{\gamma}} \right]ds \\
 &\lesssim \sup_{0< \tau \leq t}\left[\tau^{\beta} \left\|\theta(\tau) \right\|_{L^{r}}\right] \left\|\theta \right\|_{K_{\gamma,t}}\cdot  t^{-\frac{\alpha}{\gamma}},
 \end{split}
\end{equation}
where we use $\beta+\frac{\alpha}{\gamma}<1$ to apply Lemma \ref{lem:3.1} to bound the time integration by $t^{-\frac{\alpha}{\gamma}}$. By (\ref{eq:3.10}) and (\ref{eq:3.12}), for any time interval $[0,T]$ we have
\eqn \label{eq:3.13}
\left\|\mathscr{B}_{1}(v,\theta)\right\|_{K_{\gamma,T}} \lesssim \sup_{0< \tau \leq T}\left[\tau^{\beta}\|\theta(\tau)\|_{L^{q}}\right] \|\theta\|_{K_{\gamma,T}}.
\een
Since $\mathscr{B}_{2}(v,\theta)$ has the same structure as $\mathscr{B}_{1}$, we also obtain that
\eqn \label{eq:3.14}
\|\mathscr{B}_{2}(v,\theta)\|_{K_{\gamma, T}} \lesssim \sup_{0< \tau \leq T}\left[\tau^{\beta}\|\theta(\tau)\|_{L^{r}}\right] \|\theta\|_{K_{\gamma,T}}.
\een

\subsubsection{\bf Estimation of $\mathscr{B}_{3}(v,\theta)$} 
We take the $L^{p}$ norm to $\Delta_{j}\mathscr{B}_{3}(v,\theta)$.
\begin{equation}\label{eq:3.15}
 \begin{split}
 \left\|\Delta_{j} \mathscr{B}_{3}(v,\theta)(t) \right\|_{L^{p}}  &\lesssim \int^{t}_{0} \left[2^{j} e^{-(t-s)2^{\gamma j}} \sum_{k \ge j-2} \left\|\Delta_{k}\theta(s)\right\|_{L^{p}} \left\|\Delta_{k}\theta(s)\right\|_{L^{\infty}} \right]ds\\
 &\lesssim \int^{t}_{0} \left[2^{j} e^{-(t-s)2^{\gamma j}} \sum_{k \ge j-2} \left\|\Delta_{k}\theta(s) \right\|_{L^{p}} 2^{k\frac{2}{r}} \left\|\Delta_{k}\theta(s)\right\|_{L^{r}} \right]ds\\
 & \lesssim \sup_{0< \tau \leq t}\left[\tau^{\beta}\|\theta(\tau)\|_{L^{r}}\right]\int^{t}_{0} \left[2^{j} e^{-(t-s)2^{\gamma j}} s^{-\beta}\sum_{k \ge j-2} 2^{k\frac{2}{r}}\|\Delta_{k}\theta(s)\|_{L^{p}}  \right]ds,
 \end{split}
\end{equation}
where we use (\ref{eq:2.7}) to replace $\displaystyle \left\|\Delta_{k}\theta(s)\right\|_{L^{\infty}}$ by $\displaystyle 2^{k\frac{2}{r}} \left\|\Delta_{k}\theta(s)\right\|_{L^{r}}$ at the second inequality. We multiply (\ref{eq:3.15}) by $2^{j\left(\frac{2}{p}+1-\gamma\right)}$. Then,
\begin{equation} \label{eq:3.16}
 \begin{split}
 &2^{j \left(\frac{2}{p}+1-\gamma\right)} \left\|\Delta_{j} \mathscr{B}_{3}(v,\theta)(t) \right\|_{L^{p}} \\
 &\lesssim \sup_{0< \tau \leq t}\left[\tau^{\beta}\|\theta(\tau)\|_{L^{r}}\right] \\
 & \times \int^{t}_{0} \left[2^{j \left(1+\frac{2}{r}\right)} e^{-(t-s)2^{\gamma j}} s^{-\beta}\sum_{k\ge j-2} 2^{(j-k)\left(\frac{2}{p}+1-\gamma-\frac{2}{r}\right)} 2^{k\left(\frac{2}{p}+1-\gamma\right)} \left\|\Delta_{k}\theta(s)\right\|_{L^{p}} \right]ds\\
 &\lesssim \sup_{0< \tau \leq t}\left[\tau^{\beta}\|\theta(\tau)\|_{L^{r}}\right]\\
 &\times  \int^{t}_{0} \left[(t-s)^{-1/\gamma\left(1+\frac{2}{r}\right)} s^{-\beta} \sum_{k\ge j-2} 2^{(j-k)\left(\frac{2}{p}+1-\gamma-\frac{2}{r}\right)} 2^{k\left(\frac{2}{p}+1-\gamma\right)} \|\Delta_{k}\theta(s)\|_{L^{p}} \right]ds.
\end{split}
\end{equation}
We take the $l^{q}$ norm to (\ref{eq:3.16}). Since $\frac{2}{p}+1-\gamma-\frac{2}{r}>0$, by applying Young's inequality to 
\[
\sum_{k\ge j-2}a_{k-j}b_{k},  \ \text{where $a_{j}=2^{-j\left(\frac{2}{p}+1-\gamma-\frac{2}{r}\right)}$ and $b_{j}=2^{j\left(\frac{2}{p}+1-\gamma\right)} \left\|\Delta_{k}\theta(s)\right\|_{L^{p}}$},
\]
we have
\begin{equation} \label{eq:3.17}
 \begin{split}
 & \left\|\mathscr{B}_{3}(v,\theta)(t) \right\|_{\dot{B}^{\frac{2}{p}+1-\gamma}_{p,q}} \\
 & \lesssim \sup_{0< \tau \leq t}\left[\tau^{\beta}\|\theta(\tau)\|_{L^{r}}\right] \sup_{0<\tau<t}\left[\|\theta(\tau)\|_{\dot{B}^{\frac{2}{p}+1-\gamma}_{p,q}}\right]\int^{t}_{0} \left[(t-s)^{-1/\gamma(1+\frac{2}{r})} s^{-\beta} \right]ds \\
 & \lesssim \sup_{0< \tau \leq t} \left[\tau^{\beta}\|\theta(\tau)\|_{L^{r}}\right] \|\theta\|_{K_{\gamma, t}} \int^{t}_{0} \left[(t-s)^{-1/\gamma(1+\frac{2}{r})} s^{-\beta} \right]ds \lesssim \sup_{0< \tau \leq t}\left[\tau^{\beta}\|\theta(\tau)\|_{L^{r}}\right] \|\theta\|_{K_{\gamma, t}}.
 \end{split}
\end{equation}
We next estimate $\mathscr{B}_{3}(v,\theta)$ in $\dot{B}^{\frac{2}{p}+1-\gamma+\alpha}_{p,q}$. By multiplying (\ref{eq:3.15}) by $2^{j\left(\frac{2}{p}+1-\gamma+\alpha\right)}$ and following the calculation in (\ref{eq:3.16}), we have
\begin{equation} \label{eq:3.18}
 \begin{split}
 & 2^{j\left(\frac{2}{p}+1-\gamma +\alpha\right)} \left\|\Delta_{j} \mathscr{B}_{3}(v,\theta)(t)\right\|_{L^{p}} \\
 & \lesssim \sup_{0< \tau \leq t}\left[\tau^{\beta}\|\theta(\tau)\|_{L^{r}}\right] \\
 &\times \int^{t}_{0} \left[ (t-s)^{-1/\gamma\left(1+\frac{2}{r}\right)} s^{-\beta} \sum_{k\ge j-2} 2^{(j-k)\left(\frac{2}{p}+1-\gamma-\frac{2}{r}+\alpha\right)} 2^{k\left(\frac{2}{p}+1-\gamma+\alpha\right)} \left\|\Delta_{k}\theta(s)\right\|_{L^{p}}\right] ds.
 \end{split}
\end{equation}
We take the $l^{q}$ norm to (\ref{eq:3.18}). Since $\frac{2}{p}+1-\gamma-\frac{2}{r}+\alpha>0$, by applying Young's inequality to 
\[
\sum_{k\ge j-2}a_{k-j}b_{k}, \ \ \text{where $a_{j}=2^{-j\left(\frac{2}{p}+1-\gamma-\frac{2}{r}+\alpha\right)}$ and $b_{j}=2^{j\left(\frac{2}{p}+1-\gamma+\alpha\right)} \|\Delta_{k}\theta(s)\|_{L^{p}}$},
\]
we have
\begin{equation} \label{eq:3.19}
 \begin{split}
 & \left\|\mathscr{B}_{3}(v,\theta)(t)\right\|_{\dot{B}^{\frac{2}{p}+1-\gamma+\alpha}_{p,q}} \\
 & \lesssim \sup_{0< \tau \leq t}\left[\tau^{\beta}\|\theta(\tau)\|_{L^{r}}\right] \sup_{0<\tau<t}\left[\tau^{\frac{\alpha}{\gamma}} \|\theta(\tau)\|_{\dot{B}^{\frac{2}{p}+1-\gamma+\alpha}_{p,q}}\right]\int^{t}_{0} \left[ (t-s)^{-1/\gamma(1+\frac{2}{r})} s^{-\beta} s^{-\frac{\alpha}{\gamma}} \right]ds \\
 & \lesssim \sup_{0< \tau \leq t}\left[\tau^{\beta}\|\theta(\tau)\|_{L^{r}}\right] \|\theta\|_{K_{\gamma,t}} \int^{t}_{0} \left[ (t-s)^{-1/\gamma(1+\frac{2}{r})} s^{-\beta} s^{-\frac{\alpha}{\gamma}} \right]ds\\
 &\lesssim \sup_{0< \tau \leq t}\left[\tau^{\beta}\|\theta(\tau)\|_{L^{r}}\right] \|\theta\|_{K_{\gamma,t}} \cdot t^{-\frac{\alpha}{\gamma}}.
 \end{split}
\end{equation}
By (\ref{eq:3.17}) and (\ref{eq:3.19}), we conclude that for any time interval $[0,T]$
\eqn \label{eq:3.20}
\left\|\mathscr{B}_{3}(v,\theta)\right\|_{K_{\gamma, T}} \lesssim \sup_{0< \tau \leq T}\left[\tau^{\beta}\|\theta(\tau)\|_{L^{r}}\right] \|\theta\|_{K_{\gamma, T}}.
\een
In sum, by (\ref{eq:3.13}), (\ref{eq:3.14}), and (\ref{eq:3.20}), we finally have
\eqn \label{eq:3.21}
\left\|\mathscr{B}(v,\theta)\right\|_{K_{\gamma, T}} \lesssim \sup_{0< \tau \leq T}\left[\tau^{\beta}\|\theta(\tau)\|_{L^{r}}\right] \|\theta\|_{K_{\gamma, T}}.
\een
This completes the proof.

\subsection{Proof of Theorem \ref{thm:1.2} and Theorem \ref{thm:1.3}} \label{sec:3.2}
The proof of Theorem \ref{thm:1.2} and Theorem \ref{thm:1.3}, requires a couple of  elementary bounded operators which are summarized in the following two lemmas.

\begin{lemma} \label{lem:3.2}
Consider the operator 
\[
E:=e^{-[(t-s)^{1/\gamma}+s^{1/\gamma}- t^{1/\gamma}]\Lambda_{1}}
\]
for $0\leq s \leq t$.
Then $E$ is either the identity operator or has an  $L^1$ kernel whose $L^1$ norm is bounded independent of $s, t$.
\end{lemma}
\noindent
{\bf Proof.}
Clearly, 
\[
a:= (t-s)^{1/\gamma}+s^{1/\gamma}- t^{1/\gamma}
\]
is non-negative for $s \leq t$. In case $a=0$, $E$ is the identity operator, while if $a>0$, $E= e^{-a\Lambda_{1}}$ is a Fourier multiplier with symbol 
\[
\widehat{E}(\xi) = \prod_{i=1}^d e^{-a|\xi_i|}.
\]
Thus, the kernel of $E$ is given by the product of one dimensional Poisson kernels 
\[
\prod_{i=1}^d\frac{a}{\pi(a^2+x_i^2)}.
\]
The $L^1$ norm of this kernel is bounded by a constant independent of $a$.

\begin{lemma}  \label{lem:3.3}
The operator 
\[
E=e^{a^{1/\gamma}\Lambda_{1}-\frac{1}{2}a\Lambda^{\gamma}}
\]
with $\gamma>1$, is a Fourier multiplier which maps boundedly $L^p \mapsto L^p,  \ 1 < p < \infty$, and its operator norm is uniformly bounded with respect to $a \ge 0$.
\end{lemma}
\noindent
{\bf Proof.} When $a=0$, $E$ is the identity operator. When $a>0$, then $E$ is  Fourier multiplier with symbol 
\[
\widehat{E}(\xi)=e^{a^{1/\gamma}|\xi|_{1}-\frac{1}{2}a|\xi|^{\gamma}}.
\]
Since $\widehat{E}(\xi)$ is uniformly bounded for all $\xi$ and decays exponentially for $|\xi|\gg1$, the claim follows from Hormander's multiplier theorem, e.g., \cite{Stein}.

\subsubsection{\bf Proof of Theorem \ref{thm:1.2}}
We are now ready to prove a large time existence of $\Theta$ in the $L^{p}$ space. We first recall the equation of $(\Theta,V)$;
\begin{equation} \label{eq:3.22}
 \begin{split}
 \Theta(t) &=e^{t^{1/\gamma}\Lambda_{1}-t\Lambda^{\gamma}}\theta_{0}-\int^{t}_{0}\left[ e^{t^{1/\gamma}\Lambda_{1}-(t-s)\Lambda^{\gamma} } \nabla \cdot \left(e^{-s^{1/\gamma}\Lambda_{1}}V e^{-s^{1/\gamma}\Lambda_{1}} \Theta \right)(s) \right]ds\\
 & :=e^{t^{1/\gamma}\Lambda_{1}-t\Lambda^{\gamma}}\theta_{0} -\mathscr{B}(\Theta,V).
 \end{split}
\end{equation}
By Lemma \ref{lem:3.3}, it is easy to show that the linear part is equivalent to
\[
e^{-\frac{1}{2}t\Lambda^{\gamma}} \theta_{0}.
\]
As one can see at the end of the proof, the dependence on $\theta_{0}$ for the linear part is crucial to show the global existence of $\Theta$. If we show that the nonlinear part $\mathscr{B}(\Theta, V)$ is  equivalent to
\eqn \label{eq:3.23}
\int^{t}_{0}\left[ e^{-\frac{1}{2}(t-s)\Lambda^{\gamma}} \nabla \cdot (V \Theta)(s) \right]ds,
\een
we can follow \cite{Carrillo} to complete the proof. We rewrite $\mathscr{B}(V,\Theta)$ as follows.
\[
\mathscr{B}(V,\Theta)(t) = \int^{t}_{0}\left[ e^{(t^{1/\gamma}-s^{1/\gamma})\Lambda_{1}-\frac{1}{2}(t-s)\Lambda^{\gamma} } \nabla e^{-\frac{1}{2}(t-s)\Lambda^{\gamma}}  \cdot e^{s^{1/\gamma}\Lambda_{1}}\left(e^{-s^{1/\gamma}\Lambda_{1}}V e^{-s^{1/\gamma}\Lambda_{1}} \Theta\right)(s) \right]ds.
\]
We now express $(t^{1/\gamma}-s^{1/\gamma})$ as
\[
-\left((t-s)^{1/\gamma} -t^{1/\gamma} +s^{1/\gamma}\right) + (t-s)^{1/\gamma}.
\]
By Lemma \ref{lem:3.2}, $\mathscr{B}(V,\Theta)$ can be estimated by
\begin{equation}\label{eq:3.25}
 \begin{split}
 &\left\|\mathscr{B}(V,\Theta)(t)\right\|_{L^{r}} \\
 &\lesssim \int^{t}_{0} \left\|\left[ e^{(t-s)^{1/\gamma}\Lambda_{1}-\frac{1}{2}(t-s)\Lambda^{\gamma}} \nabla e^{-\frac{1}{2}(t-s)\Lambda^{\gamma}}  \cdot e^{s^{1/\gamma}\Lambda_{1}} \left(e^{-s^{1/\gamma}\Lambda_{1}}V e^{-s^{1/\gamma}\Lambda_{1}} \Theta\right)(s) \right]\right\|_{L^{r}}ds.
 \end{split}
\end{equation}
By Lemma \ref{lem:3.3}, the right-hand side of (\ref{eq:3.25}) can be replaced by
\eqn \label{eq:3.26}
\left\|\mathscr{B}(V,\Theta)(t)\right\|_{L^{r}} \lesssim \int^{t}_{0} \left\|\left[ \nabla e^{-\frac{1}{2}(t-s)\Lambda^{\gamma}} \cdot e^{s^{1/\gamma}\Lambda_{1}} \left(e^{-s^{1/\gamma}\Lambda_{1}}V e^{-s^{1/\gamma}\Lambda_{1}} \Theta \right)(s) \right]\right\|_{L^{r}}ds.
\een
Using the fact that
\[
\left\|\nabla e^{-t\Lambda^{\gamma}}f \right\|_{L^{r}} \lesssim t^{-1/\gamma-\frac{2}{\gamma}\left(\frac{1}{q}-\frac{1}{r}\right)} \left\|f\right\|_{L^{q}},
\]
we estimate (\ref{eq:3.26}) as
\eqn \label{eq:3.27}
\left\|\mathscr{B}(V,\Theta)(t)\right\|_{L^{r}} \lesssim \int^{t}_{0} (t-s)^{-1/\gamma-\frac{2}{\gamma}\left(\frac{1}{q}-\frac{1}{r}\right)} \left\|e^{s^{1/\gamma}\Lambda_{1}}\left(e^{-s^{1/\gamma}\Lambda_{1}}V e^{-s^{1/\gamma}\Lambda_{1}} \Theta\right)(s) \right\|_{L^{q}}ds.
\een
To estimate the right-hand side of (\ref{eq:3.27}), we introduce the bilinear operators $B_{t}$ of the form
\begin{equation}\label{eq:3.28}
 \begin{split}
 B_{t}(f,g) &:=e^{t^{1/\gamma}\Lambda_{1}} \left(e^{-t^{1/\gamma}\Lambda_{1}} f e^{-t^{1/\gamma}\Lambda_{1}} g\right)\\
 &=\int_{R^{2}} \int_{R^{2}} e^{ix\cdot(\xi+\eta)} e^{t^{1/\gamma}(|\xi+\eta|_{1}-|\xi|_{1}-|\eta|_{1})} \hat{f}(\xi) \hat{g}(\eta)d\xi d\eta.
 \end{split}
\end{equation}
Recall that for a vector $\xi=(\xi_1, \xi_2)$, we denoted $|\xi|_1=|\xi_1|+|\xi_{2}|$. For $\xi=(\xi_1,\xi_2), \eta = (\eta_1, \eta_2)$, we split the domain of integration of the above integral into sub-domains depending on the sign of $\xi_j, \eta_j$ and $\xi_j+\eta_j$. In order to do so, we introduce the operators acting on one variable  (see page 253 in \cite{Rieusset}) by
\[
K_1f := \frac{1}{2\pi} \int_0^{\infty} e^{\imath x \xi}\hat{f}(\xi)\, d\xi, \quad K_{-1}f := \frac{1}{2\pi} \int_{-\infty}^{0} e^{\imath x \xi}\hat{f}(\xi)\, d\xi.
\]
Let the operators $L_{t,-1}$ and $L_{t,1}$ be defined by
\[
L_{t,1}f = f, \quad L_{t,-1}f = \frac{1}{2\pi}\int_{R} e^{\imath x\xi}e^{-2t|\xi|}\hat{f}(\xi)\, d\xi.
\]
For $ \vec{\alpha}=(\alpha_1, \alpha_2)$, $\vec{\beta}  =(\beta_1, \beta_2) \in \{-1,1\}^2$, denote the operator
\[
Z_{t,\vec \alpha , \vec \beta} = K_{\beta_1}L_{t,\alpha_1 \beta_1} \otimes \cdots \otimes K_{\beta_2}L_{t,\alpha_2\beta_2}, \quad K_{\vec \alpha} = k_{\alpha_1} \otimes K_{\alpha_2}.
\]
The above tensor product  means that the $j-$th operator in the tensor product acts on the $j-$th variable of the function $f(x_1, x_2)$.
A tedious (but elementary) calculation now yields the following identity:
\eqn \label{eq:3.29}
B_{t}(f,g)=\sum_{(\vec \alpha, \vec \beta, \vec \gamma) \in \{-1,1\}^{2\times 2}} K_{\alpha_1} K_{\alpha_2} \left( Z_{t, \vec \alpha, \vec \beta}f Z_{t, \vec \alpha, \vec \gamma}g \right).
\een
We note that the operators $K_{\vec \alpha}$, $Z_{t, \vec \alpha , \vec \beta}$ defined above, being linear combinations of Fourier multipliers (including Hilbert transform)  and the identity operator, commute with $\Lambda_{1}$ and $\Lambda$. Moreover, they are bounded linear operators on $L^p, 1 < p < \infty$ and the corresponding operator norm of $Z_{t,\vec \alpha, \vec \beta}$ is bounded independent of $t\ge 0$. Therefore,
\eqn \label{eq:3.30}
\left\|B_{t}(f,g)\right\|_{L^{q}} \lesssim \left\|fg\right\|_{L^{q}}.
\een
We apply the above argument to the right-hand side of (\ref{eq:3.27}) to conclude that
\eqn \label{eq:3.31}
\left\|\mathscr{B}(V,\Theta)(t)\right\|_{L^{r}} \lesssim \int^{t}_{0} (t-s)^{-1/\gamma-\frac{2}{\gamma}(\frac{1}{q}-\frac{1}{r})} \left\|V(s)\Theta(s)\right\|_{L^{q}}ds.
\een
We now follow the proof in \cite{Carrillo} line by line to obtain 
\eqn
\left\|\Theta\right\|_{G_{\gamma, T}} \lesssim \|\theta_{0}\|_{L^{\frac{2}{\gamma-1}}} + \sup_{0< \tau \leq T}\left[\tau^{\beta}\|\Theta(\tau)\|_{L^{r}}\right] \|\Theta\|_{G_{\gamma, T}}.
\een
This a priori estimate implies a local existence of a solution in $G_{\gamma,T}$ on $[0,T_{0}]$ for some $T_{0}>0$. Then, we restart the problem with initial data at $T_{0}/2$. Since 
\[
\left\|\theta\left(T_{0}/2\right) \right\|_{L^{\frac{2}{\gamma-1}}} \leq \|\theta_{0}\|_{L^{\frac{2}{\gamma-1}}}
\]
by the maximum principle (\cite{Cordoba}), we have the solution on $\left[T_{0}/2, 3T_{0}/2\right]$. Repeating this process, we can reach any pre-assigned Time $T$ in finitely many steps.

\subsubsection{\bf Proof of Theorem \ref{thm:1.3}}
As for the proof of Theorem \ref{thm:1.1}, we only need to obtain the following a priori estimate:
\eqn \label{eq:3.32}
\left\|\Theta\right\|_{K_{\gamma, T}} \lesssim \|\theta_{0}\|_{\dot{B}^{\frac{2}{p}+1-\gamma}_{p,q}} + \sup_{0< \tau \leq T}\left[\tau^{\beta}\|\Theta(\tau)\|_{L^{r}}\right] \left\|\Theta\right\|_{K_{\gamma, T}}.
\een
By Lemma \ref{lem:3.3}, the linear estimation is obvious. By repeating the argument in the proof of Theorem \ref{thm:1.2}, we can estimate the nonlinear term as of $\theta$. For the reader's convenience, we provide a few lines. We take $\Delta_{j}$ to $\mathscr{B}(V,\Theta)$ in (\ref{eq:3.22}) and take the $L^{p}$ norm. By Lemma \ref{lem:3.2} and \ref{lem:3.3},
\begin{equation}  \label{eq:3.33}
 \begin{split}
 \left\| \Delta_{j}e^{t^{1/\gamma}\Lambda_{1}}\mathscr{B}(V, \Theta) \right\|_{L^{p}}  \lesssim \int^{t}_{0} \left[ e^{-\frac{1}{2}(t-s)2^{\gamma j}} 2^{j} \left\|e^{s^{1/\gamma}\Lambda_{1}} \Delta_{j} \big( e^{-s^{1/\gamma}\Lambda_{1}}V e^{-s^{1/\gamma}\Lambda_{1}}\Theta)(s)\big) \right\|_{L^{p}} \right]ds.
 \end{split}
\end{equation}
We decompose the product $e^{-s^{1/\gamma}\Lambda_{1}}V e^{-s^{1/\gamma}\Lambda_{1}}\Theta$ as paraproduct:
\[
T_{\left(e^{-s^{1/\gamma}\Lambda_{1}}V\right)} e^{-s^{1/\gamma}\Lambda_{1}}\Theta +T_{\left(e^{-s^{1/\gamma}\Lambda_{1}}\Theta\right)} e^{-s^{1/\gamma}\Lambda_{1}}V +R\left(e^{-s^{1/\gamma}\Lambda_{1}}V, e^{-s^{1/\gamma}\Lambda_{1}}\Theta \right).
\]
Then,
\begin{equation}  \label{eq:3.34}
 \begin{split}
 & \left\| \Delta_{j} e^{t^{1/\gamma}\Lambda_{1}}\mathscr{B}(V,\Theta) \right\|_{L^{p}}  \\
 & \lesssim \int^{t}_{0}  \left[ e^{-\frac{1}{2}(t-s)2^{\gamma j}} 2^{j} \left\|e^{s^{1/\gamma}\Lambda_{1}} \left( e^{-s^{1/\gamma}\Lambda_{1}}S_{j}V  e^{-s^{1/\gamma}\Lambda_{1}}\Delta_{j}\Theta)(s)\right) \right\|_{L^{p}} \right]ds\\
 & + \int^{t}_{0}  \left[ e^{-\frac{1}{2}(t-s)2^{\gamma j}} 2^{j} \left\|e^{s^{1/\gamma}\Lambda_{1}} \left( e^{-s^{1/\gamma}\Lambda_{1}}S_{j}\Theta  e^{-s^{1/\gamma}\Lambda_{1}}\Delta_{j}V)(s)\right) \right\|_{L^{p}} \right]ds\\
 & + \int^{t}_{0} \sum_{k \ge j-2}\left[ e^{-\frac{1}{2}(t-s)2^{\gamma j}} 2^{j} \left\| e^{s^{1/\gamma}\Lambda_{1}} \left( e^{-s^{1/\gamma}\Lambda_{1}}\Delta_{k}V  e^{-s^{1/\gamma}\Lambda_{1}}\Delta_{k}\Theta)(s)\right) \right\|_{L^{p}} \right]ds
 \end{split}
\end{equation}
By using (\ref{eq:3.30}), we can estimate (\ref{eq:3.34}) as
\begin{equation}  \label{eq:3.35}
 \begin{split}
 & \left\| \Delta_{j} e^{t^{1/\gamma}\Lambda_{1}}\mathscr{B}(V,\Theta) \right\|_{L^{p}} \\
 & \lesssim  \int^{t}_{0}  \left[ e^{-\frac{1}{2}(t-s)2^{\gamma j}} 2^{j} \left( \left\|S_{j}V(s) \Delta_{j}\Theta(s) \right\|_{L^{p}} + \left\| S_{j}\Theta(s) \Delta_{j}V(s) \right\|_{L^{p}} \right)\right]ds\\
 & + \int^{t}_{0} e^{-\frac{1}{2}(t-s)2^{\gamma j}} 2^{j} \sum_{k \ge j-2}\left[\left\| \Delta_{k}V(s) \Delta_{k}\Theta(s) \right\|_{L^{p}} \right]ds
 \end{split}
\end{equation}
Therefore, we can follow the calculations line by line from (\ref{eq:3.6}) to (\ref{eq:3.21}) in the proof of Theorem \ref{thm:1.1} to complete the proof.

\section{Critical case: Proof of Theorem \ref{thm:1.4} and Theorem \ref{thm:1.5}} \label{sec:4}
To prove Theorem \ref{thm:1.4} and Theorem \ref{thm:1.5}, we need several lemmas. First, we need the following representation to estimate $v$ in $L^{\infty}$.

\begin{lemma} \label{lem:4.1}
Oseen Kernel \cite{Rieusset}: The operator $O_{t}=\mathscr{R}K_{t}$ is a convolution operator whose kernel $\tilde{K}_{t}$ satisfies 
\[
\tilde{K}_{t}(x)=\frac{1}{t^{d}}\tilde{K}\left(x/t\right)
\]
for a smooth function $\tilde{K}$ such that for all $\alpha\in \mathbb{N}^{d}$, 
\[
(1+|x|)^{d+\alpha}\partial^{\alpha}\tilde{K}\in L^{\infty}.
\]
In particular, $\tilde{K}$ is in $L^{p}$ for $p>1$.
\end{lemma}

Next, we provide the $L^{p}$ bounds of the Poisson kernel and its Oseen kernel in two dimensions. The proof is easily obtained by their representation.
\begin{lemma} \label{lem:4.2}
For any $1< p < \infty$,
\[
\left\|P_{t}\right\|_{L^{p}}\lesssim t^{-2\left(1-\frac{1}{p}\right)}, \quad \left\|\mathscr{R}P_{t}\right\|_{L^{p}}\lesssim t^{-2\left(1-\frac{1}{p}\right)}.
\]
\end{lemma}

To deal with the time singularities appearing in Lemma \ref{lem:4.2}, we need the following lemma.
\begin{lemma} \label{lem:4.3}
Hardy-Littlewood-Sobolev Inequality \cite{Lieb}: Let $0<\lambda <d$, $\frac{1}{p} +\frac{\lambda}{d} +\frac{1}{q} =2$. Then,
\[
\left|\int_{\mathbb{R}^{d}} \int_{\mathbb{R}^{d}} \frac{f(x)g(y)}{|x-y|^{\lambda}} dydx\right| \lesssim \|f\|_{L^{p}} \|g\|_{L^{q}}.
\]
In particular, for the one dimensional case,
\[
\sup_{t>0}\left|\int^{t}_{0} \frac{1}{|t-s|^{\frac{1}{2}}}a(s)ds\right| \lesssim \|a\|_{L^{2}}.
\]
\end{lemma}

\subsection{Existence: Proof of Theorem \ref{thm:1.4}} 
For the reader's convenience, we will repeat the computation in \cite{Bae2}.

\subsubsection{\bf $\dot{H}^{1}$ bound} 
By taking one derivative $\nabla$ to (\ref{eq:1.1}),
\eqn \label{eq:4.1}
\nabla \theta_{t}+v\cdot \nabla \nabla \theta +\nabla v \cdot \nabla \theta-\Lambda_{1} \nabla \theta=0.
\een
We multiply (\ref{eq:4.1}) by $\nabla \theta$ and integrate over $\mathbb{R}^{2}$. Then,
\eqn \label{eq:4.2}
\frac{1}{2}\frac{d}{dt} \left\|\nabla \theta\right\|^{2}_{L^{2}}+\left\|\nabla^{\frac{3}{2}} \theta \right\|^{2}_{L^{2}} \leq \int_{\mathbb{R}^{2}} \left|\nabla v\right| \left|\nabla \theta \right| \left|\nabla v \right|dx \leq \left\|\nabla \theta\right\|_{L^{2}} \left\|\nabla v\right\|^{2}_{L^{4}}.
\een
By the Sobolev embedding $\dot{H}^{1/2} \subset L^{4}$ in two dimensions, we can replace $\left\|\nabla v\right\|_{L^{4}}$ in (\ref{eq:4.2}) by $\left\|\nabla^{\frac{3}{2}}v\right\|_{L^{2}}$. Then,
\[
\frac{d}{dt} \left\|\nabla \theta\right\|^{2}_{L^{2}}+\left\|\nabla^{\frac{3}{2}}\theta \right\|^{2}_{L^{2}}\lesssim \left\|\nabla \theta\right\|_{L^{2}} \left\|\nabla^{\frac{3}{2}} v\right\|^{2}_{L^{2}}.
\]
Since $\displaystyle v=\left(-\mathscr{R}_{2}\theta, \mathscr{R}_{1}\theta\right)$, and the Riesz transforms  are bounded in $H^{s}$, we have
\eqn \label{eq:4.3}
\frac{d}{dt} \left\|\nabla \theta \right\|^{2}_{L^{2}}+\left\|\nabla^{\frac{3}{2}}\theta \right\|^{2}_{L^{2}}\lesssim \left\|\nabla \theta \right\|_{L^{2}} \left\|\nabla^{\frac{3}{2}} \theta \right\|^{2}_{L^{2}}.
\een
Integrating (\ref{eq:4.3}) in time, we have
\eqn \label{eq:4.4}
\left\|\nabla \theta \right\|^{2}_{L^{\infty}_{t}L^{2}}+\left\|\nabla^{\frac{3}{2}}\theta \right\|^{2}_{L^{2}_{t}L^{2}} \lesssim \left\|\nabla\theta_{0} \right\|^{2}_{L^{2}} +\left\|\nabla \theta \right\|_{L^{\infty}_{t}L^{2}} \left\|\nabla^{\frac{3}{2}}\theta \right\|^{2}_{L^{2}_{t}L^{2}}.
\een

\subsubsection{\bf $L^{\infty}$ bound} 
To obtain the $L^{\infty}$ bound of $\theta$, we express $\theta$ as the integral form:
\eqn  \label{eq:4.5}
\theta(t)=e^{-t\Lambda}\theta_{0} - \int^{t}_{0} \left[\nabla e^{-(t-s)\Lambda} \cdot \left(v \theta\right)(s) \right]ds.
\een
By taking the $L^{\infty}$ norm, we have
\begin{equation} \label{eq:4.6}
 \begin{split}
 \left\|\theta(t) \right\|_{L^{\infty}} & \leq \left\|\theta_{0}\right\|_{L^{\infty}}+\int^{t}_{0} \left\|e^{-(t-s)\Lambda}(v\cdot\nabla \theta)(s) \right\|_{L^{\infty}}ds \\
&\lesssim \left\|\theta_{0}\right\|_{L^{\infty}}+\int^{t}_{0} \left\|e^{-(t-s)\Lambda} \right\|_{L^{\frac{4}{3}}} \left\|v(s)\right\|_{L^{\infty}} \left\|\nabla \theta(s) \right\|_{L^{4}}ds.
 \end{split}
\end{equation}
By Lemma \ref{lem:4.2},
\eqn \label{eq:4.7}
\left\|\theta(t)\right\|_{L^{\infty}} \lesssim \left\|\theta_{0}\right\|_{L^{\infty}} +\int^{t}_{0}\frac{1}{\sqrt{t-s}} \left\|v(s)\right\|_{L^{\infty}} \left\|\nabla \theta(s)\right\|_{L^{4}}ds.
\een
By Lemma \ref{lem:4.3} and the Sobolev embedding $\dot{H}^{1/2}\subset L^{4}$ in two dimensions, we finally have
\begin{equation} \label{eq:4.8}
 \begin{split}
 \left\|\theta(t) \right\|_{L^{\infty}}  &\leq \left\|\theta_{0}\right\|_{L^{\infty}}+\left\|v\right\|_{L^{\infty}_{t}L^{\infty}} \left\|\nabla \theta\right\|_{L^{2}_{t}L^{4}} \lesssim \left\|\theta_{0}\right\|_{L^{\infty}}+\left\|v\right\|_{L^{\infty}_{t}L^{\infty}} \left\|\nabla^{\frac{3}{2}} \theta \right\|_{L^{2}_{t}L^{2}}.
 \end{split}
\end{equation}
We next estimate  $v$ in $L^{\infty}$. For simplicity, we set $v=\mathscr{R}\theta$ and take $\mathscr{R}$ to (\ref{eq:4.5}).
\[
v(t)=e^{-t\Lambda} v_{0}-\int^{t}_{0} \left[\mathscr{R} e^{-(t-s)\Lambda}(v\cdot\nabla \theta)(s)\right]ds.
\]
By following the estimates of $\theta$ obtained in (\ref{eq:4.6}) to (\ref{eq:4.8}), we have
\begin{equation} \label{eq:4.9}
 \begin{split}
 \left\|v(t) \right\|_{L^{\infty}} &\lesssim \|v_{0}\|_{L^{\infty}}+\int^{t}_{0} \left\|\mathscr{R}e^{-(t-s)\Lambda} \right\|_{L^{\frac{4}{3}}} \left\|v(s)\right\|_{L^{\infty}} \left\|\nabla \theta(s)\right\|_{L^{4}}ds\\
 & \lesssim \left\|v_{0} \right\|_{L^{\infty}}+\int^{t}_{0}\frac{1}{\sqrt{t-s}} \left\|v(s)\right\|_{L^{\infty}} \left\|\nabla \theta(s)\right\|_{L^{4}}ds\\
 & \lesssim \left\|v_{0}\right\|_{L^{\infty}}+\left\|v\right\|_{L^{\infty}_{t}L^{\infty}} \left\|\nabla^{\frac{3}{2}} \theta\right\|_{L^{2}_{t}L^{2}}.
 \end{split}
\end{equation}
In sum, by (\ref{eq:4.4}), (\ref{eq:4.8}), and (\ref{eq:4.9}), we obtain that
\begin{equation} \label{eq:4.10}
 \begin{split}
 & \left\|\theta\right\|_{L^{\infty}_{t}L^{\infty}}+\left\|v\right\|_{L^{\infty}_{t}L^{\infty}}+\left\|\theta\right\|_{L^{\infty}_{t}\dot{H}^{1}}+\left\|\nabla^{\frac{3}{2}}\theta\right\|_{L^{2}_{t}L^{2}}\\
&\lesssim \left\|\theta_{0}\right\|_{L^{\infty}}+\left\|v_{0}\right\|_{L^{\infty}}+\left\|\theta_{0}\right\|_{\dot{H}^{1}} + \left(\left\|\theta\right\|_{L^{\infty}_{t}L^{\infty}}+\left\|v\right\|_{L^{\infty}_{t}L^{\infty}}+\left\|\theta\right\|_{L^{\infty}_{t}\dot{H}^{1}}+\left\|\nabla^{\frac{3}{2}}\theta\right\|_{L^{2}_{t}L^{2}}\right)^{2},
 \end{split}
\end{equation}
which implies that the existence of a global-in-time solution $\theta$ and $v$ in $L^{\infty}_{t}(L^{\infty}\cap \dot{H}^{1})\cap {L}^{2}_{t}\dot{H}^{\frac{3}{2}}$ provided that $\|\theta_{0}\|_{L^{\infty}}+\|v_{0}\|_{L^{\infty}}+\|\theta_{0}\|_{\dot{H}^{1}}$ is sufficiently small.

\subsubsection{\bf $\dot{B}^{\frac{2}{p}}_{p, q}$ bound} 
We take $\Delta_{j}$ to (\ref{eq:4.5}) and take the $L^{p}$ norm. By (\ref{eq:2.8}), we have
\eqn \label{eq:4.11}
\left\|\Delta_{j}\theta(t)\right\|_{L^{p}} \lesssim e^{-t2^{j}} \left\|\Delta_{j}\theta_{0}\right\|_{L^{p}} +\int^{t}_{0} \left[e^{-(t-s)2^{j}} \left\|\Delta_{j}(v \cdot \nabla \theta)(s)\right\|_{L^{p}}\right]ds.
\een
By taking the $L^{\infty}$ norm in time, we obtain
\eqn \label{eq:4.12}
\left\|\Delta_{j}\theta\right\|_{L^{\infty}_{t}L^{p}} \lesssim \left\|\Delta_{j}\theta_{0}\right\|_{L^{p}} +\left\|\Delta_{j}(v \cdot \nabla \theta)\right\|_{L^{1}_{t}L^{p}}.
\een
We multiply (\ref{eq:4.12}) by $2^{j\frac{2}{p}}$ and take the $l^{q}$ norm. Then,
\eqn \label{eq:4.13}
\left\|\theta\right\|_{\tilde{L}^{\infty}_{t}\dot{B}^{\frac{2}{p}}_{p,q}}\lesssim \left\|\theta_{0}\right\|_{\dot{B}^{\frac{2}{p}}_{p,q}} + \left\|v\cdot\nabla \theta\right\|_{\tilde{L}^{1}_{t}\dot{B}^{\frac{2}{p}}_{p,q}}.
\een
Similarly, by taking the $L^{1}$ norm in time to (\ref{eq:4.12}), multiplying by $2^{j\left(\frac{2}{p}+1\right)}$, and taking the $l^{q}$ norm, we obtain
\eqn
\label{eq:4.14}
\left\|\theta\right\|_{\tilde{L}^{1}_{t}\dot{B}^{\frac{2}{p}+1}_{p,q}}\lesssim \left\|\theta_{0}\right\|_{\dot{B}^{\frac{2}{p}}_{p,q}} + \left\|v\cdot\nabla \theta\right\|_{\tilde{L}^{1}_{t}\dot{B}^{\frac{2}{p}}_{p,q}}.
\een
By adding (\ref{eq:4.13}) and (\ref{eq:4.14}), we finally have
\eqn \label{eq:4.15}
\|\theta\|_{\tilde{L}^{\infty}_{t}\dot{B}^{\frac{2}{p}}_{p,q}} + \|\theta\|_{\tilde{L}^{1}_{t} \dot{B}^{\frac{2}{p}+1}_{p,q}}\lesssim \|\theta_{0}\|_{\dot{B}^{\frac{2}{p}}_{p,q}} + \left\|v\cdot\nabla \theta\right\|_{\tilde{L}^{1}_{t}\dot{B}^{\frac{2}{p}}_{p,q}}.
\een
We now estimate the nonlinear term 
\[
\left\|v\cdot\nabla \theta\right\|_{\tilde{L}^{1}_{t}\dot{B}^{\frac{2}{p}}_{p,q}}
\]
by using the paraproduct of $v$ and $\theta$. By (\ref{eq:2.6}),
\eqn \label{eq:4.17}
\Delta_{j}(v\theta)=S_{j}v\Delta_{j}\theta +S_{j}\theta \Delta_{j}v +\sum_{k\ge j-2}\Delta_{k}v\Delta_{k}\theta.
\een
By taking the $L^{p}$ norm to (\ref{eq:4.17}), we have
\begin{equation} \label{eq:4.18}
 \begin{split}
 \left\|\Delta_{j}(v\theta) \right\|_{L^{p}} &\lesssim \left\|S_{j}v\right\|_{L^{\infty}} \left\|\Delta_{j}\theta\right\|_{L^{p}} + \left\|S_{j}\theta\right\|_{L^{\infty}} \left\|\Delta_{j}v\right\|_{L^{p}} + \sum_{k\ge j-2} \left\|\Delta_{k}v\right\|_{L^{\infty}} \left\|\Delta_{k}\theta\right\|_{L^{p}} \\
 & =\left\|S_{j}v\right\|_{L^{\infty}} \left\|\Delta_{j}\theta \right\|_{L^{p}} + \left\|S_{j}\theta\right\|_{L^{\infty}} \left\|\Delta_{j}v\right\|_{L^{p}} \\
 &+ \sum_{k\ge j-2}2^{-k\left(\frac{2}{p}+1\right)} \left\|\Delta_{k}v\right\|_{L^{\infty}} 2^{k\left(\frac{2}{p}+1\right)} \left\|\Delta_{k}\theta\right\|_{L^{p}}.
 \end{split}
\end{equation}
By taking the $L^{1}$ norm in time to (\ref{eq:4.18}),
\begin{equation} \label{eq:4.19}
 \begin{split}
 \left\|\Delta_{j}(v\theta)\right\|_{L^{1}_{t}L^{p}} &\lesssim \left\|S_{j}v\right\|_{L^{\infty}_{t}L^{\infty}} \left\|\Delta_{j}\theta\right\|_{L^{1}_{t}L^{p}} + \left\|S_{j}\theta\right\|_{L^{\infty}_{t}L^{\infty}} \left\|\Delta_{j}v\right\|_{L^{1}_{t}L^{p}} \\
&+ \sum_{k\ge j-2}2^{-k\left(\frac{2}{p}+1\right)} \left\|\Delta_{k}v\right\|_{L^{\infty}_{t}L^{\infty}} 2^{k\left(\frac{2}{p}+1\right)} \left\|\Delta_{k}\theta\right\|_{L^{1}_{t}L^{\infty}}.
 \end{split}
\end{equation}
We multiply (\ref{eq:4.19}) by $2^{j\left(\frac{2}{p}+1\right)}$. Then,
\begin{equation} \label{eq:4.20}
 \begin{split}
 2^{j\left(\frac{2}{p}+1\right)} \left\|\Delta_{j}(v\theta)\right\|_{L^{1}_{t}L^{p}} &\lesssim 2^{j\left(\frac{2}{p}+1\right)} \left( \left\|S_{j}v\right\|_{L^{\infty}_{t}L^{\infty}} \left\|\Delta_{j}\theta\right\|_{L^{1}_{t}L^{p}} + \left\|S_{j}\theta \right\|_{L^{\infty}_{t}L^{\infty}} \left\|\Delta_{j}v\right\|_{L^{1}_{t}L^{p}}\right) \\
&+ \sum_{k\ge j-2}2^{(j-k)\left(\frac{2}{p}+1\right)} \left\|\Delta_{k}v\right\|_{L^{\infty}_{t}L^{\infty}} 2^{k\left(\frac{2}{p}+1\right)} \left\|\Delta_{k}\theta \right\|_{L^{1}_{t}L^{\infty}}\\
&\lesssim \left\|v\right\|_{L^{\infty}_{t}L^{\infty}} 2^{j\left(\frac{2}{p}+1\right)} \left\|\Delta_{j}\theta\right\|_{L^{1}_{t}L^{p}} + \left\|\theta\right\|_{L^{\infty}_{t}L^{\infty}} 2^{j\left(\frac{2}{p}+1\right)} \left\|\Delta_{j}v\right\|_{L^{1}_{t}L^{p}} \\
&+ \left\|v\right\|_{L^{\infty}_{t}L^{\infty}} \sum_{k\ge j-2}2^{(j-k)\left(\frac{2}{p}+1\right)} 2^{k\left(\frac{2}{p}+1\right)} \left\|\Delta_{k}\theta\right\|_{L^{1}_{t}L^{\infty}}.
 \end{split}
\end{equation}
We take the $l^{q}$ norm to (\ref{eq:4.20}). Since $\frac{2}{p}+1>0$, by applying Young's inequality to 
\[
\sum_{k\ge j-2}a_{k-j}b_{k}, \ \ \text{where $a_{j}=2^{-j\left(\frac{2}{p}+1\right)}$ and $b_{j}=2^{j\left(\frac{2}{p}+1\right)} \left\|\Delta_{k}\theta\right\|_{L^{1}_{t}L^{p}}$},
\]
we have
\eqn \label{eq:4.21}
\left\|v\cdot\nabla \theta\right\|_{\tilde{L}^{1}_{t}\dot{B}^{\frac{2}{p}}_{p,q}} \lesssim \left\|v\right\|_{{L}^{\infty}_{t}L^{\infty}} \left\|\theta\right\|_{\tilde{L}^{1}_{t}\dot{B}^{\frac{2}{p}+1}_{p,q}} + \left\|\theta\right\|_{{L}^{\infty}_{t}L^{\infty}} \left\|v\right\|_{\tilde{L}^{1}_{t}\dot{B}^{\frac{2}{p}+1}_{p,q}}.
\een
Since $v=\left(-\mathscr{R}_{2}\theta, \mathscr{R}_{1}\theta\right)$, and Besov spaces are bounded under the Riesz transformations,
\eqn \label{eq:4.22}
\left\|v\cdot\nabla \theta \right\|_{\tilde{L}^{1}_{t}\dot{B}^{\frac{2}{p}}_{p,q}} \lesssim \left(\left\|\theta\right\|_{{L}^{\infty}_{t}L^{\infty}}+\left\|v\right\|_{{L}^{\infty}_{t}L^{\infty}}\right) \left\|\theta\right\|_{\tilde{L}^{1}_{t}\dot{B}^{\frac{2}{p}+1}_{p,q}}.
\een
Combining (\ref{eq:4.15}) and (\ref{eq:4.22}), we finally have
\eqn \label{eq:4.23}
\left\|\theta\right\|_{\tilde{L}^{\infty}_{t}\dot{B}^{\frac{2}{p}}_{p,q}} + \|\theta\|_{\tilde{L}^{1}_{t}\dot{B}^{\frac{2}{p}+1}_{p,q}}\lesssim \|\theta_{0}\|_{\dot{B}^{\frac{2}{p}}_{p,q}} +\left(\|\theta\|_{{L}^{\infty}_{t}L^{\infty}} + \|v\|_{{L}^{\infty}_{t}L^{\infty}}\right)\|\theta\|_{\tilde{L}^{1}_{t}\dot{B}^{\frac{2}{p}+1}_{p,q}},
\een
which completes the a priori estimate of the solution in $\tilde{L}^{\infty}_{t}\dot{B}^{\frac{2}{p}}_{p,q} \cap \tilde{L}^{1}_{t}\dot{B}^{\frac{2}{p}+1}_{p,q}$ by using the smallness of $ \|\theta\|_{{L}^{\infty}_{t}L^{\infty}} + \|v\|_{{L}^{\infty}_{t}L^{\infty}}$ derived in (\ref{eq:4.10}).

\subsection{Analyticity: Proof of Theorem \ref{thm:1.5}} 
We first recall
\eqn \label{eq:4.24}
\Theta(t)=e^{\frac{1}{4}t\Lambda_{1}}\theta(t), \quad V(t)=e^{\frac{1}{4}t\Lambda_{1}}v(t)
\een
and the equation they satisfy:
\eqn \label{eq:4.25}
\Theta(t)=e^{\frac{1}{4}t\Lambda_{1}-t\Lambda}\theta_{0}-\int^{t}_{0}\left[ e^{\frac{1}{4}t\Lambda_{1}-(t-s)\Lambda} \left(e^{-\frac{1}{4}s\Lambda_{1}}V \cdot \nabla e^{-\frac{1}{4}s\Lambda_{1}} \Theta\right)(s) \right]ds.
\een
In order to reduce (\ref{eq:4.25}) to
\eqn \label{eq:4.26}
\Theta(t)=e^{-\frac{1}{2}t\Lambda_{1}}\theta_{0}-\int^{t}_{0}\left[ e^{-\frac{1}{2}(t-s)\Lambda_{1}} \left(V \cdot \nabla \Theta\right)(s) \right]ds
\een
in terms of estimation, we need to show that the variation of Lemma \ref{lem:3.2} and Lemma \ref{lem:3.3} for $\gamma =1$ works. Lemma \ref{lem:3.2} is identically applicable to the case $\gamma=1$. Moreover, 
\[
e^{\frac{1}{4}t\Lambda_{1}-\frac{1}{2}t\Lambda}
\]
is a Fourier multiplier which maps boundedly $L^p \mapsto L^p,  \ 1 < p < \infty$, and its operator norm is uniformly bounded with respect to $t \ge 0$. Therefore, along the lines of the estimation of $\theta$ in  $L^{\infty}$ and $\dot{B}^{\frac{2}{p}}_{p,q}$, we can find the following two bounds immediately from (\ref{eq:4.26}):
\begin{equation} \label{eq:4.27}
 \begin{split}
 & \left\|\Theta\right\|_{L^{\infty}_{t}L^{\infty}}+\left\|V\right\|_{L^{\infty}_{t}L^{\infty}} \lesssim \|\theta_{0}\|_{L^{\infty}}+ \|v_{0}\|_{L^{\infty}}+ \left(\|\Theta\|_{L^{\infty}_{t}L^{\infty}}+\|V\|_{L^{\infty}_{t}L^{\infty}} \right) \left\|\nabla^{\frac{3}{2}}\Theta \right\|_{L^{2}_{t}L^{2}},\\
 & \left\|\Theta\right\|_{\tilde{L}^{\infty}_{t}\dot{B}^{\frac{2}{p}}_{p,q}} + \left\|\Theta\right\|_{\tilde{L}^{1}_{t}\dot{B}^{\frac{2}{p}+1}_{p,q}}\lesssim \|\theta_{0}\|_{\dot{B}^{\frac{2}{p}}_{p,q}} +\left(\|\Theta\|_{{L}^{\infty}_{t}L^{\infty}} + \|V\|_{{L}^{\infty}_{t}L^{\infty}}\right)\|\Theta\|_{\tilde{L}^{1}_{t}\dot{B}^{\frac{2}{p}+1}_{p,q}}.
  \end{split}
\end{equation}
It remains to derive the $\dot{H}^{1}$ estimation of $\Theta$, which will be obtained by the energy method.
\begin{equation} \label{eq:4.29}
 \begin{split}
 & \frac{1}{2}\frac{d}{dt}\int_{R^{2}} \left|\nabla\Theta \right|^{2}dx = \int_{R^{2}} \left[\nabla \Theta \cdot \nabla\Theta_{t}\right]dx=\int_{R^{2}} \left[\nabla\Theta \cdot\nabla\left(e^{\frac{1}{4}t\Lambda_{1}}\theta\right)_{t}\right]dx\\
 &= \int_{R^{2}} \left[\nabla\Theta \cdot\nabla\left(\frac{1}{4}\Lambda_{1}\Theta -e^{\frac{1}{4}t\Lambda_{1}} (v\cdot\nabla \theta)-\Lambda \Theta\right) \right]dx \\
 &=\int_{R^{2}} \left[\nabla\Theta \cdot \nabla \left(\frac{1}{4}\Lambda_{1}-\Lambda \right)\Theta\right] dx -\int_{R^{2}} \left[\nabla\Theta \cdot \nabla e^{\frac{1}{4}t\Lambda_{1}} \left(e^{-\frac{1}{4}t\Lambda_{1}}V \cdot \nabla e^{-\frac{1}{4}t\Lambda_{1}} \Theta\right)\right]dx.
  \end{split}
\end{equation}
By using $\frac{1}{4}|\xi|_{1}<\frac{1}{2}|\xi|$ and the boundedness of $B_{t}(V,\Theta)$,
\eqn \label{eq:4.30}
\frac{1}{2}\frac{d}{dt}\int_{R^{2}}|\nabla\Theta|^{2}dx \lesssim -\frac{1}{2}\|\nabla^{\frac{3}{2}} \Theta\|^{2}_{L^{2}} +\|\nabla^{\frac{3}{2}} \Theta\|_{L^{2}} \|\nabla^{\frac{3}{2}}\left(\Theta \mathscr{R}\Theta \right)\|_{L^{2}}.
\een
By applying the product rule (\cite{Kato}) to $\Lambda_{1}^{\frac{3}{2}}\left(\Theta \mathscr{R}\Theta \right)$, we finally have
\eqn \label{eq:4.31}
\frac{1}{2}\frac{d}{dt}\int_{R^{2}} \left|\nabla\Theta\right|^{2}dx \lesssim  -\frac{1}{2} \left\|\nabla^{\frac{3}{2}} \Theta\right\|^{2}_{L^{2}} +\left\|\nabla^{\frac{3}{2}} \Theta \right\|^{2}_{L^{2}} \left(\|\Theta\|_{L^{\infty}} +\|V\|_{L^{\infty}}\right).
\een
Integrating (\ref{eq:4.31}) in time,
\eqn \label{eq:4.32}
\left\|\nabla \Theta \right\|^{2}_{L^{\infty}_{t}L^{2}} +\left\|\nabla^{\frac{3}{2}} \Theta \right\|^{2}_{L^{2}_{t}L^{2}} \lesssim \left\|\nabla\theta_{0}\right\|^{2}_{L^{2}} + \left(\|\Theta\|_{L^{\infty}_{t}L^{\infty}} +\|V\|_{L^{\infty}_{t}L^{\infty}} \right) \left\|\nabla^{\frac{3}{2}} \Theta\right\|^{2}_{L^{2}_{t}L^{2}}.
\een
By (\ref{eq:4.27}) and (\ref{eq:4.32}), we complete the proof.

\section*{Acknowledgments}

H.B. gratefully acknowledges the support by the Center for Scientific Computation and Mathematical Modeling (CSCAMM) at University of Maryland where this research was performed. A. B. also thanks the Center for Scientific Computation and Mathematical Modeling at the University of Maryland for hosting him while this work was completed and the University of North Carolina at Charlotte for partial support via the reassignment of duties program. Research was supported in part by NSF grant DMS10-08397 and ONR grant 000141210318 (H. Bae and E. Tadmor) and DMS 11-09532 (A. Biswas).

\bibliographystyle{amsplain}

\end{document}